\documentclass[12pt,a4paper]{article}
\usepackage{color}
\usepackage{amsmath, amssymb, latexsym}
\usepackage{amsthm}
\usepackage{graphicx}
\usepackage{epsf}
\usepackage{fullpage}
\usepackage[pdftex,pagebackref]{hyperref}

\newtheorem{theorem}{Theorem}[section]
\newtheorem*{theorem*}{Theorem}
\newtheorem{lemma}[theorem]{Lemma}
\newtheorem{proposition}[theorem]{Proposition}
\newtheorem{assumption}[theorem]{Assumption}

\newtheorem{property}[theorem]{Property}

\newtheorem{remark}[theorem]{Remark}

\allowdisplaybreaks
\newcommand{\beq}{\begin{equation}}
\newcommand{\eeq}{\end{equation}}

\def \cL {\mathcal L}

\newcommand {\ar}{\rightarrow}

\numberwithin{equation}{section}

\newcommand{\noib}{\noindent $\bullet $~}
\newcommand{\ie}{i.e.\;}

\newcommand{\resp}{\emph{resp. }}

\newcommand{\piq}{\frac{\pi}{4}}

\newcommand{\hH}{\widehat{H}}

\def \R {{\mathbb R}}

\newcommand{\bS}{\mathbb{S}}

\newcommand{\cE}{\mathcal{E}}
\newcommand{\cF}{\mathcal{F}}

\newcommand{\cS}{\mathcal{S}}

\newcommand{\lambdah}{\hat{\lambda}}

\title{Some nodal properties of the quantum harmonic oscillator and other Schr\"{o}dinger operators in $\R^2$}

\author{P. B\'erard \\ Institut Fourier, Universit\'{e} Grenoble Alpes and CNRS, B.P.74,
\\ F38402 Saint Martin d'H\`{e}res Cedex, France.\\
\href{mailto:pierrehberard@gmail.com}{pierrehberard@gmail.com}\\
and \\
B. Helffer\\
Laboratoire Jean Leray, Universit\'{e} de Nantes and CNRS,\\
F44332 Nantes Cedex France, and \\
Laboratoire de Math\'ematiques, Univ. Paris-Sud 11.\\
\href{mailto:Bernard.Helffer@univ-nantes.fr}{Bernard.Helffer@univ-nantes.fr}}

\date{\small June 21, 2015, revised July 31, 2016}

\begin{document}
\maketitle

\begin{abstract} {For the spherical Laplacian on the sphere and
for the Dirichlet Laplacian in the square}, Antonie
Stern claimed in her  PhD thesis (1924) the existence of an
infinite sequence of eigenvalues whose corresponding eigenspaces
contain an eigenfunction with exactly two nodal domains. These  results were  given complete proofs respectively by Hans Lewy in
1977, and the authors in 2014 (see also Gauthier-Shalom--Przybytkowski, 2006). In this paper, we obtain similar results for the two dimensional isotropic quantum harmonic oscillator. In the opposite direction, we construct an infinite sequence of regular eigenfunctions with as many nodal domains as allowed by Courant's theorem, up to a factor $\frac{1}{4}$. A classical question for a $2$-dimensional bounded domain is to estimate the length of the nodal set of a Dirichlet eigenfunction in terms of the square root of the energy. In the last section, we consider some Schr\"{o}dinger operators $-\Delta + V$ in $\R^2$ and we provide bounds for the length of the nodal set of an eigenfunction with energy $\lambda$ in the classically permitted region $\{V(x) < \lambda\}$.
\end{abstract}


\section{Introduction and main results}\label{S-intro}

Given a finite interval $]a,b[$ and a continuous function $q:[a,b]\mapsto \R$, consider the one-dimensional self-adjoint eigenvalue problem
\begin{equation}\label{E-sturm}
-y'' + qy = \lambda y \text{~in~} ]a,b[\,, ~~ y(a) = y(b)=0\,.
\end{equation}
Arrange the eigenvalues in increasing order, $\lambda_1(q) < \lambda_2(q) < \cdots$. A classical theorem of C.~Sturm \cite{Hin05} states that an eigenfunction $u$ of \eqref{E-sturm} associated with $\lambda_k(q)$ has exactly $(k-1)$ zeros in $]a,b[$ or, equivalently, that the zeros of $u$ divide $]a,b[$ into $k$ sub-intervals.

In higher dimensions, one can consider the eigenvalue problem for the Laplace-Beltrami operator $-\Delta_g$ on a compact connected Riemannian manifold $(M,g)$, with Dirichlet condition in case $M$ has a boundary $\partial M$,
\begin{equation}\label{E-higher}
-\Delta u = \lambda u \text{~in~} M, ~~ u|_{\partial M} = 0.
\end{equation}
Arrange the eigenvalues in non-decreasing order, with multiplicities,
\begin{equation*}
\lambda_1(M,g) < \lambda_2(M,g) \le \lambda_3(M,g) \le \ldots
\end{equation*}
Denote by $M_0$ the interior of $M$, $M_0 := M \setminus \partial M$.

Given an eigenfunction $u$ of $-\Delta_g$, denote by
\begin{equation}\label{S1-2}
N(u) := \overline{ \left\lbrace x \in M_0 ~|~ u(x)=0 \right\rbrace} \end{equation}
the \emph{nodal set} of $u$, and by
\begin{equation}\label{S1-4}
\mu(u) := \# \left\lbrace \text{~connected components of~} M_0 \setminus N(u) \right\rbrace
\end{equation}
the number of \emph{nodal domains} of $u$ i.e., the number of connected components of the complement of $N(u)$.

Courant's theorem \cite{Cou} states that if $-\Delta_g u =\lambda_k(M,g) u$, then $\mu(u) \le k$.\medskip

In this paper, we investigate three natural questions about Courant's theorem in the framework of the 2D isotropic quantum harmonic oscillator.\medskip

\textbf{Question~1}.~ In view of Sturm's theorem, it is natural to ask whether Courant's upper bound is sharp, and to look for lower bounds for the number of nodal domains, depending on the geometry of $(M,g)$ and the eigenvalue. Note that for orthogonality reasons, for any $k \ge 2$ and any eigenfunction associated with $\lambda_k(M,g)$, we have $\mu(u) \ge 2$.

We shall say that $\lambda_k(M,g)$ is \emph{Courant-sharp} if there exists an eigenfunction $u$, such that $-\Delta_g u = \lambda_k(M,g) u$ and $\mu(u) = k$. Clearly, $\lambda_1(M,g)$ and $\lambda_2(M,g)$ are always Courant-sharp eigenvalues. Note that if $\lambda_3(M,g) = \lambda_2(M,g)$, then $\lambda_3(M,g)$ is not Courant-sharp.

The first results concerning Question~1 were stated by Antonie Stern in her 1924 PhD thesis \cite{St}  written  under the supervision of R.~Courant.

\newpage
\begin{theorem}\label{T-stern}[A.~Stern, {\cite{St}}]
\vspace{-3mm}
\begin{enumerate}
   \item For the square $[0,\pi]\times[0,\pi]$ with Dirichlet boundary condition, there is a sequence of eigenfunctions $\{u_r, r\ge 1\}$ such that
       $$
       -\Delta u_r = (1+4r^2) u_r\,, \text{~and~} \mu(u_r)=2\,.
       $$
   \item For the sphere $\bS^2$, there exists a sequence of eigenfunctions $u_{\ell}, \ell \ge 1$ such that
       $$
       -\Delta_{\bS^2} u_{\ell} = \ell (\ell + 1)u_{\ell}\,,
       \text{~and~} \mu(u_{\ell}) = 2 \text{~or~} 3\,,
       $$
       depending on whether $\ell$ is odd or even.
 \end{enumerate}
\end{theorem}%

Stern's arguments are not fully satisfactory. In 1977, H. Lewy \cite{Lew} gave a complete independent proof for the case of the sphere, without any reference to \cite{St} (see also \cite{BeHe2}). More recently, the authors \cite{BeHe} gave a complete proof for the case of the square with Dirichlet boundary conditions (see also Gauthier-Shalom--Przybytkowski \cite{GSP}).

The original motivation of this paper was to investigate the possibility to extend Stern's results to the case of the two-dimensional isotropic quantum harmonic oscillator  $\hH := - \Delta + |x|^2$ acting on $L^2(\R^2,\R)$ (we will say ``harmonic oscillator'' for short). After the publication of the first version of this paper \cite{BeHe-I}, T.~Hoffmann-Ostenhof informed us of the unpublished master degree thesis of J.~Leydold \cite{LeyD}.

An orthogonal basis of eigenfunctions of the harmonic oscillator $\hH$ is given by
\begin{equation}
\phi_{m,n} (x,y) = H_m(x) H_n (y) \exp(- \frac{x^2+y^2}{2})\,,
\end{equation}
for $(m,n)\in \mathbb N^2$, where $H_n$ denotes the
Hermite polynomial of degree $n$. For Hermite polynomials, we use the definitions and notation of Szeg\"o \cite[\S 5.5]{Sz}.

The eigenfunction $\phi_{m,n}$ corresponds to the eigenvalue
$2(m +n+1)\,,$
\begin{equation}
\hH \phi_{m,n} = 2(m+n+1)\, \phi_{m,n}\, .
\end{equation}
It follows that the eigenspace $\cE_{n}$ of $\widehat{H}$ associated with the eigenvalue $\lambdah(n)= 2(n +1)$ has dimension $(n + 1)$, and is generated by the eigenfunctions $\phi_{n,0}$, $\phi_{n-1,1}$,
\dots, $\phi_{0,n}$.\medskip

We summarize Leydold's main results in the following theorem.

\begin{theorem}\label{T-Leydold}[J.~Leydold, \cite{LeyD}]
\vspace{-3mm}
\begin{enumerate}
  \item For $n \ge 2$, and for any nonzero $u \in \cE_n$,
  \begin{equation}\label{Ley1}
  \mu(u) \le \mu_n^L :=\frac{n^2}{2} + 2\,.
  \end{equation}
  \item The lower bound on the number of nodal domains is given by
  \begin{equation}\label{Ley2}
  \min \left\lbrace \mu(u) ~|~ u \in \cE_n, u \neq 0 \right\rbrace =
  \left\{
    \begin{array}{ll}
    1 & \text{if~} n = 0\,,\\
    3 & \text{if~} n \equiv 0 \pmod{4}\,, n\ge 4\,,\\
    2 & \text{if~} n \not \equiv 0 \pmod{4}\,.
    \end{array}
  \right.
  \end{equation}
\end{enumerate}
\end{theorem}%

\textbf{Remark.} When $n\ge 3$, the estimate \eqref{Ley1} is better than Courant's bound which is $ \mu_n^C := \frac{n^2}{2} + \frac{n}{2} + 1$. The idea of the proof is to apply Courant's method separately to odd and to even eigenfunctions (with respect to the map $x \mapsto - x$). A consequence of \eqref{Ley1} is that the only Courant-sharp eigenvalues of the harmonic oscillator are the first, second and fourth eigenvalues.  The same ideas work for the sphere as well \cite{LeyD, Ley}. For the analysis of Courant-sharp eigenvalues of the square with Dirichlet boundary conditions, see \cite{Pl,BeHe}.\medskip

Leydold's proof that there exist eigenfunctions of the harmonic oscillator satisfying \eqref{Ley2} is quite involved. In this paper, we give a simple proof that, for any odd integer $n$, there exists a one-parameter family of eigenfunctions with exactly two nodal domains in $\cE_n$. More precisely, for $\theta \in [0,\pi[$, we consider the following curve in $\cE_n$,
\begin{equation}
\Phi^{\theta}_{n} := \cos\theta \, \phi_{n,0} + \sin\theta \,
\phi_{0,n}\,,
\end{equation}
\begin{equation*}
\Phi^{\theta}_{n}(x,y) = \left( \cos\theta \, H_n(x) + \sin\theta \, H_n(y) \right) \, \exp(- \frac{x^2+y^2}{2})\,.
\end{equation*}

We prove the following theorems (Sections~\ref{S-th1} and \ref{S-th2}).

\begin{theorem}\label{th1}   Assume that  $n$ is odd. Then,
there exists an open interval $I_{\frac{\pi}{4}}$ containing
$\frac{\pi}{4}$, and an open interval $I_{\frac{3\pi}{4}}$,
containing $\frac{3\pi}{4}$, such that for
$$
\theta \in I_{\frac{\pi}{4}} \cup I_{\frac{3\pi}{4}} \setminus
\{\frac{\pi}{4},\frac{3\pi}{4}\}\,,
$$
the nodal set $N(\Phi^{\theta}_n)$ is a connected simple regular
curve, and the eigenfunction $\Phi^{\theta}_n$ has two nodal domains
in $\R^2$.
\end{theorem}

\begin{theorem}\label{th2}
 Assume that $n$ is odd. Then, there exists $\theta_c >0$ such
that, for $0 < \theta < \theta_c\,$, the nodal set
$N(\Phi^{\theta}_n)$ is a connected simple regular curve, and the
eigenfunction $\Phi^{\theta}_n$ has two nodal domains in $\R^2$.
\end{theorem}

\textbf{Remark}. The value $\theta_c$ and the intervals can be computed numerically. The proofs of the theorems actually show that $]0,\theta_c[\, \cap I_{\frac{\pi}{4}} = \emptyset$.\medskip

\textbf{Question~2}. How good/bad is Courant's upper bound on the number of nodal domains? Consider the eigenvalue problem \eqref{E-higher}. Given $k \ge 1$, define $\mu(k)$ to be the maximum value of $\mu(u)$ when $u$ is an eigenfunction associated with the eigenvalue $\lambda_k(M,g)$. Then,
\begin{equation}\label{E-pleijel}
\limsup_{k \mapsto \infty} \frac{\mu(k)}{k} \le \gamma(n)
\end{equation}
where $\gamma(n)$ is a universal constant which only depends on the dimension $n$ of $M$. Furthermore, for $n \ge 2$, $\gamma(n) < 1$. The idea of the proof, introduced by Pleijel in 1956, is to use a Faber-Krahn type isoperimetric inequality and Weyl's asymptotic law. Note that the constant $\gamma(n)$ is not sharp. For more details and references, see \cite{Pl, Do14}.

As a corollary of the above result,  the eigenvalue problem \eqref{E-higher} has only finitely many Courant-sharp eigenvalues.

The above result  gives a quantitative  improvement of Courant's theorem in the case of the Dirichlet Laplacian in a bounded open set  of $\R^2$. When trying to implement the strategy of Pleijel for the harmonic oscillator, we get into trouble because of the absence of a reasonable Faber-Krahn inequality. P. Charron \cite{Cha,Cha15} has obtained the following theorem.

\begin{theorem}
If $(\lambda_n,u_n)$ is an infinite sequence of eigenpairs of $\widehat H$, then
\begin{equation}
\lim\sup \frac{\mu (u_n)}{n} \leq \gamma(2) = \frac{4}{j_{0,1}^2}\,,
\end{equation}
where $ j_{0,1}$ is the first positive zero of the Bessel function  of order $0\,$.
\end{theorem}

This is in some sense surprising that the statement is exactly the same as in the case of  the Dirichlet realization of the Laplacian in a bounded open set in $\mathbb R^2$. The proof  does actually not use the isotropy of the harmonic potential, can be extended to the $n$-dimensional case, but strongly  uses the explicit knowledge of the eigenfunctions. We refer to \cite{ChHeHO16} for further results in this direction.\medskip

A related question concerning the estimate \eqref{E-pleijel} is whether the order of magnitude is correct. In the case of the $2$-sphere, using spherical coordinates, one can find decomposed spherical harmonics $u_{\ell}$ of degree $\ell$, with associated eigenvalue $\ell (\ell + 1)$, such that $\mu(u_{\ell}) \sim \frac{\ell^2}{2}$ when $\ell$ is large, whereas Courant's upper bound is equivalent to $\ell^2$. These spherical harmonics have critical zeros and their nodal sets have self-intersections. In \cite[Section~2]{EJN}, the authors construct spherical harmonics $v_{\ell}$, of degree $\ell$, without critical zeros i.e., whose nodal sets are disjoint closed regular curves, such that $\mu(v_{\ell}) \sim \frac{\ell^2}{4}$. These spherical harmonics $v_{\ell}$ have as many nodal domains as allowed by Courant's theorem, up to a factor $\frac 14$. Since the $v_{\ell}$'s are regular, this property is stable under small perturbations in the same eigenspace.

In Section~\ref{S-ovals}, in a direction opposite to Theorems~\ref{th1} and \ref{th2}, we construct eigenfunctions of the harmonic oscillator $\widehat{H}$ with ``many'' nodal domains.

\begin{theorem}\label{T-mnd}
For the harmonic oscillator $\widehat{H}$ in $L^2(\R^2)$, there exists a sequence of eigenfunctions $\{u_k, k\ge 1\}$ such that $\widehat{H}u_k = \hat{\lambda}(k) u_k$, with $\hat{\lambda}(4k) = 2 (4k+1)$, $u_k$ as no critical zeros (i.e. has a regular nodal set), and
$$
\mu(u_k) \sim \frac{(4k)^2}{8}.
$$
\end{theorem}%

\textbf{Remarks}.\\
1.~The above estimate is, up to a factor $\frac{1}{4}$, asymptotically the same as the upper bounds for the number of nodal domains given by Courant and Leydold.\\
2.~A related question is to analyze the zero set when $\theta$ is a random variable. We refer to \cite{HZZ} for  results in this direction. The above questions are related to the question of spectral minimal partitions \cite{HHOT}. In the case of the harmonic oscillator similar questions appear in the analysis of the properties of
ultracold atoms (see for example \cite{RL}).\medskip

\textbf{Question~3}. Consider the eigenvalue problem \eqref{E-higher}, and assume for simplicity that $M$ is a bounded domain in $\R^2$. Fix any small number $r$, and a point $x \in M$ such that $B(x,r) \subset M$. Let $u$ be a Dirichlet eigenfunction associated with the eigenvalue $\lambda$ and assume that $\lambda \ge \frac{\pi j_{0,1}^2}{r^2}$. Then $N(u) \cap B(x,r) \not = \emptyset$. This fact follows from the monotonicity of the Dirichlet eigenvalues, and indicates that the length of the nodal set should tend to infinity as the eigenvalue tends to infinity. The first results in this direction are due to Br\"{u}ning and Gromes \cite{BrGr, Br78} who show that the length of the nodal set $N(u)$ is bounded from below by a constant times $\sqrt{\lambda}$.  For further results in this direction (Yau's conjecture), we refer to \cite{DF, LJ, LoMa16, SZ1, SZ2}.
In Section~\ref{S-lbnz}, we investigate this question for the harmonic oscillator.

\begin{theorem}\label{lbnz-P1}
Let $\delta \in ]0,1[$ be given. Then, there exists a positive constant $C_{\delta}$ such that for $\lambda$ large enough, and for any nonzero eigenfunction of the isotropic $2D$ quantum harmonic oscillator,
\begin{equation}\label{lbnz-2}
\widehat{H} := - \Delta + |x|^2\,, ~~ \widehat{H} u = \lambda u\,,
\end{equation}
the length of $N(u) \cap B\left( \sqrt{\delta \lambda} \right)$ is bounded from below by $C_{\delta}\, \lambda^{\frac{3}{2}}$.
\end{theorem}

As a matter of fact, we prove a lower bound for more general Schr\"{o}dinger operators in $\R^2$ (Propositions~\ref{lbnz-P2} and \ref{prop7.9}), shedding some light on the exponent $\frac{3}{2}$ in the above estimate. In Section~\ref{SS-lj}, we investigate upper and lower bounds for the length of the nodal sets, using the method of Long Jin \cite{LJ}.

{\bf Acknowledgements.} The authors would like to thank P.~Charron for providing an earlier copy of \cite{Cha}, T.~Hoffmann-Ostenhof for pointing out \cite{LeyD}, J.~Leydold for providing a copy of his master degree thesis \cite{LeyD}, Long Jin for enlightening discussions concerning \cite{LJ} and Section~\ref{S-lbnz}, I.~Polterovich for pointing out \cite{Cha}, as well as D.~Jakobson and M.~Persson-Sundqvist for useful comments.  During the preparation of this paper, B.~H. was Simons foundation visiting fellow at the Isaac Newton Institute in Cambridge.
Finally, the authors would like to thank the referee for his remarks which helped improve the paper.

\section{A  reminder on Hermite polynomials}\label{S-hermite}

We use the definition, normalization, and notation of Szeg\"{o}'s book
\cite{Sz}. With these choices, the Hermite polynomial $H_n$ has the following properties,
\cite[\S~5.5 and Theorem~6.32]{Sz}.\vspace{-3mm}

\begin{enumerate}
\item  $H_n$ satisfies the differential equation
$$
y''(t) - 2 t \,y'(t) + 2n\, y(t) = 0\,.
$$
\item  $H_n$ is a polynomial of degree $n$ which is even (\resp odd) for $n$
even (\resp odd).
\item  $H_n(t) = 2 t\, H_{n-1}(t) - 2 (n-1) \,
H_{n-2}(t)\,, ~n\ge 2\,, ~~H_0(t)= 1\,, ~~H_1(t) = 2t\,.$
\item $H_n$ has $n$ simple zeros $t_{n,1} < t_{n,2} < \cdots <
t_{n,n}\,$.
\item
$$
 H_n(t) =  2 t\, H_{n-1}(t) - H'_{n-1}(t)\,.
$$
\item \begin{equation}
H'_n(t) = 2n H_{n-1} (t)\,.
\end{equation}
\item The coefficient of $t^n$ in  $H_n$ is $2^n$.
\item
$$
\int _{-\infty}^{+\infty} e^{-t^2} |H_n(t)| ^2\, dt = \pi^\frac 12 \,2^n\, n! \,.
$$
\item The  first zero $t_{n,1}$ of $H_n$ satisfies
\begin{equation}
 t_{n,1} = (2n+1)^\frac 12 - 6^{-\frac 12} (2n+1)^{-\frac 16}
(i_1+\epsilon_n)\,,
\end{equation}
where $i_1$ is the first positive real zero of the Airy function,
and $\lim_{n\rightarrow +\infty} \epsilon_n =0\,$.
\end{enumerate}

The following result (Theorem 7.6.1 in Szeg\"o's book \cite{Sz}) will
also be useful.

\begin{lemma}\label{L-leh}
The successive relative maxima of $t \mapsto |H_n(t)|$ form an increasing
sequence for $t \geq 0$\,.
\end{lemma}
{\bf Proof.}
\\
It is enough to observe that the function
$$
\Theta_n (t):= 2 n H_n(t)^2 + H_n'(t)^2
$$
satisfies
$$
\Theta_n'(t) = 4 t \, (H'_n(t)) ^2 \,.
$$\hfill $\square$

\section{Stern-like constructions for the harmonic oscillator in the case $n$-odd}\label{S-sl-odd}

\subsection{The case of the square}

Consider the square  $[0,\pi]^2$, with Dirichlet boundary
conditions, and the following families of eigenfunctions associated
with the eigenvalues $\hat{\lambda}(1,2r) := 1+ 4r^2$, where $r$ is
a positive integer, and $\theta \in [0,\pi/4]$,
$$
(x,y) \mapsto \cos \theta \, \sin x \, \sin(2ry)  + \sin \theta \,
\sin(2rx) \, \sin y\,.
$$

\begin{figure}[!ht]
\begin{center}
\includegraphics[width=15cm]{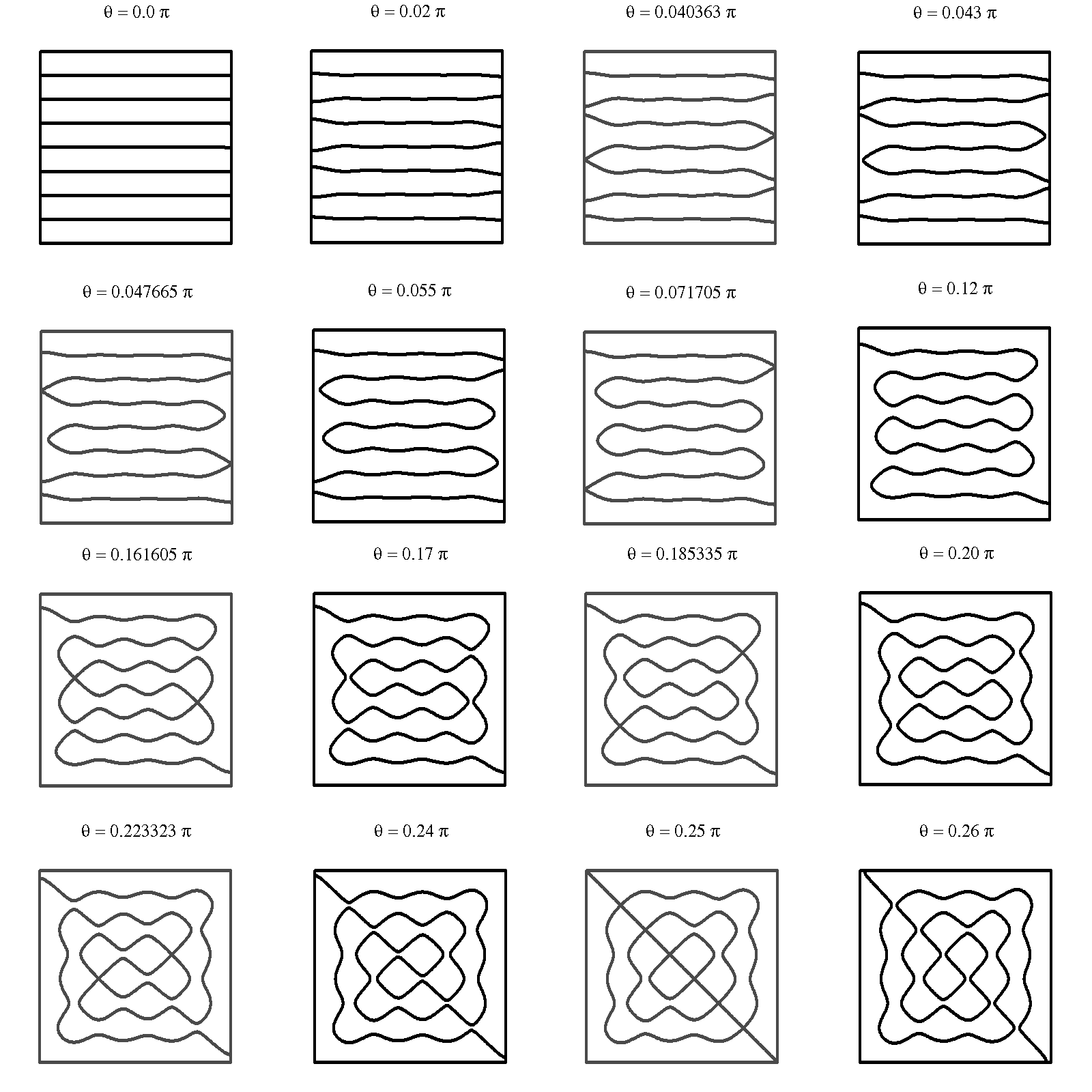}
\caption{Nodal sets for the Dirichlet eigenvalue $\hat{\lambda}(1,8)$ of the square.} \label{sq-1}
\end{center}
\end{figure}

According to \cite{St}, for any given $r \ge 1$, the typical evolution of
the nodal sets when $\theta$ varies is similar to the case $r=4$
shown in Figure~\ref{sq-1} \cite[Figure~6.9]{BeHe}. Generally
speaking, the nodal sets deform continuously, except for \emph{finitely many} values of $\theta$, for which self-intersections of the nodal set appear or disappear or, equivalently, for which critical zeros of the eigenfunction appear/disappear.

We would like to get similar results for the isotropic quantum
harmonic oscillator.

\subsection{Symmetries}

 Recall the notation,
\begin{equation}\label{P2}
\Phi^{\theta}_n(x,y) := \cos\theta \,
\phi_{n,0}(x,y) + \sin\theta \, \phi_{0,n}(x,y)\,.
\end{equation}
 More simply,
$$
\Phi^{\theta}_n(x,y) = \exp\left( -\frac{x^2+y^2}{2} \right)\, \left( \cos\theta \, H_n(x) + \sin\theta \, H_n(y) \right)\,.
$$

Since $\Phi^{\theta + \pi}_n = - \Phi^{\theta}_n$, it suffices to
vary the parameter $\theta$ in the interval $[0,\pi[$.

\emph{Assuming $n$ is odd}, we have the following symmetries.

\begin{equation}\label{P4}
\left\lbrace
\begin{array}{ll}
\Phi^{\theta}_n(-x,y) & = \Phi^{\pi-\theta}_n(x,y)\,,\\[4pt]
\Phi^{\theta}_n(x,-y) & = - \Phi^{\pi-\theta}_n(x,y)\,,\\[4pt]
\Phi^{\theta}_n(y,x) & = \Phi^{\frac{\pi}{2}-\theta}_n(x,y)\,.\\[4pt]
\end{array}
\right.
\end{equation}

When $n$ is odd, it therefore suffices to vary the parameter
$\theta$ in the interval $[0,\piq]$. The case $\theta = 0$ is
particular, so that we shall mainly consider $\theta \in ]0,\piq]$.

\subsection{Critical zeros}

A \emph{critical zero} of $\Phi^{\theta}_n$ is a point $(x,y) \in
\R^2$ such that both $\Phi^{\theta}_n$ and its differential
$d\Phi^{\theta}_n$ vanish at $(x,y)$. The critical zeros of
$\Phi^{\theta}_n$ satisfy the following equations.

\begin{equation}\label{P6}
\left\lbrace
\begin{array}{ll}
\cos\theta \, H_n(x) + \sin\theta \, H_n(y) &= 0\,,\\[4pt]
\cos\theta \, H'_n(x) & = 0 \,,\\[4pt]
\sin\theta \, H'_n(y) &= 0\,.\\
\end{array}
\right.
\end{equation}

Equivalently, using the properties of the Hermite polynomials, a
point $(x,y)$ is a critical zero of $\Phi^{\theta}_n$ if and only if

\begin{equation}\label{P8}
\left\lbrace
\begin{array}{ll}
\cos\theta \, H_n(x) + \sin\theta \, H_n(y) &= 0\,,\\[4pt]
\cos\theta \, H_{n-1}(x) & = 0 \,,\\[4pt]
\sin\theta \, H_{n-1}(y) &= 0\,.\\
\end{array}
\right.
\end{equation}

The only  \emph{possible} critical zeros of the eigenfunction
$\Phi^{\theta}_n$ are the points $(t_{n-1,i}\,,t_{n-1,j})$ for $1
\le i,j \le (n-1)$, where the coordinates are the zeros of the
Hermite polynomial $H_{n-1}$. The point $(t_{n-1,i}\,,t_{n-1,j})$ is
a critical zero of $\Phi^{\theta}_n$ if and only if $\theta =
\theta(i,j)\,$, where $\theta(i,j) \in ]0,\pi[$ is uniquely determined by the equation,
\begin{equation}\label{P10}
\cos\left( \theta(i,j) \right) \, H_n(t_{n-1,i}) + \sin\left(
\theta(i,j) \right) \, H_n(t_{n-1,j}) = 0\,.
\end{equation}
The values $\theta(i,j)$ will be called \emph{critical values} of the parameter $\theta$, the other values \emph{regular values}. Here we have used the fact that $H_n$ and $H'_n$ have no common
zeros. We have proved the following lemma.

\begin{lemma}\label{L-cz}
For $\theta \in [0,\pi[$, the eigenfunction $\Phi^{\theta}_n$ has no
critical zero, unless $\theta$ is one of the critical values $\theta(i,j)$ defined
by equation \eqref{P10}. In particular $\Phi^{\theta}_n$ has no
critical zero, except for finitely many values of the parameter
$\theta \in [0,\pi[$. Given a pair $(i_0,j_0) \in \{1, \ldots, n-1 \}$, let $\theta_0 = \theta(i_0,j_0)$,
 be defined by \eqref{P10} for the pair $(t_{n-1,i_0}\,, t_{n-1,j_0})\,$. Then, the function
$\Phi^{\theta_0}_n$ has finitely many critical zeros, namely the
points $ (t_{n-1,i}\,, t_{n-1,j})$ which satisfy
\begin{equation}
\cos\theta_0 \, H_n(t_{n-1,i}) + \sin\theta_0 \, H_n(t_{n-1,j}) =
0\,,
\end{equation}
among them the point  $(t_{n-1,i_0}\,, t_{n-1,j_0})\,$.
\end{lemma}

\textbf{Remarks}.\\
 From the general properties of nodal lines
\cite[Properties~5.2]{BeHe}, we  derive the following facts.
\vspace{-3mm}
\begin{enumerate}
    \item When $\theta \not \in \left\lbrace \theta(i,j) ~|~ 1 \le i,j \le n-1
    \right\rbrace$, the nodal set $N(\Phi^{\theta}_n)$ of the eigenfunction $\Phi^{\theta}_n$, is a smooth $1$-dimensional submanifold of $\R^2$ (a collection of pairwise distinct connected simple regular curves).
    \item When $\ \theta  \in \left\lbrace \theta(i,j)~|~ 1 \le i,j \le n-1
    \right\rbrace$, the nodal set $N(\Phi^{\theta}_n)$ has   finitely
    many singularities which are double crossings\footnote{This result is actually general for any eigenfunction of the harmonic oscillator, as stated in \cite{LeyD}, on the basis of Euler's formula and Courant's theorem.}. Indeed, the Hessian
    of the
    function $\Phi^{\theta}_n$ at a \emph{critical zero}
    $(t_{n-1,i},t_{n-1,j})$ is
    given by
    \begin{equation*}\label{hessian}
    \mathrm{Hess}_{(t_{n-1,i},t_{n-1,j})}\Phi^{\theta}_n =
    \exp{(-\frac{t_{n-1,i}^2 + t_{n-1,j}^2}{2})}\,
    \begin{pmatrix}
      \cos\theta \, H''_n(t_{n-1,i}) & 0 \\
      0 & \sin\theta \, H''_n(t_{n-1,j}) \\
    \end{pmatrix}\,,
    \end{equation*}
and the assertion follows from the fact that $H_{n-1}$ has simple
zeros.
\end{enumerate}

\subsection{General properties of the nodal set $N(\Phi^{\theta}_n)$}\label{ss-gpns}

Denote by $\cL$ the finite lattice
\begin{equation}\label{P22}
\cL := \left\lbrace (t_{n,i}\,,t_{n,j}) ~|~ 1 \le i,j \le n
\right\rbrace \subset \R^2 \,,
\end{equation}
consisting of points whose coordinates are the zeros of the Hermite
polynomial $H_n$.

The horizontal and vertical lines $\{y=t_{n,i}\}$ and $\{x=t_{n,j}\}$, $1 \le i,j \le n$, form a checkerboard like pattern in $\R^2$ which can be colored according to the sign of the function $H_n(x)\, H_n(y)$ (grey where the function is positive, white where it is negative). We will refer to the following properties as \emph{the checkerboard argument}, compare with \cite{St, BeHe}.

For symmetry reasons, we can assume that $\theta \in ]0,\piq]$.

(i) We have the following inclusions for the nodal  sets $N(\Phi^{\theta}_n)\,$,
\begin{equation}\label{P24}
\cL \subset N(\Phi^{\theta}_n) \subset \cL \cup \left\lbrace (x,y)
\in \R^2 ~|~ H_n(x) \, H_n(y) < 0 \right\rbrace \,.
\end{equation}

(ii) The nodal set $N(\Phi^{\theta}_n)$ does not meet the vertical
lines $\{x=t_{n,i}\}$, or the horizontal lines $\{y=t_{n,i}\}$ away
from the set $\cL$.

(iii) The lattice point $(t_{n,i},t_{n,j})$ is not a
critical zero of $\Phi^{\theta}_n$ (because $H_n$ and $H'_n$
have no common zero). As a matter of fact, near a lattice point, the
nodal set $N(\Phi^{\theta}_n)$ is a single arc through the lattice
point, with a tangent which is neither horizontal, nor
vertical.\medskip

\begin{figure}[ht!]
\begin{center}
\includegraphics[width=10cm]{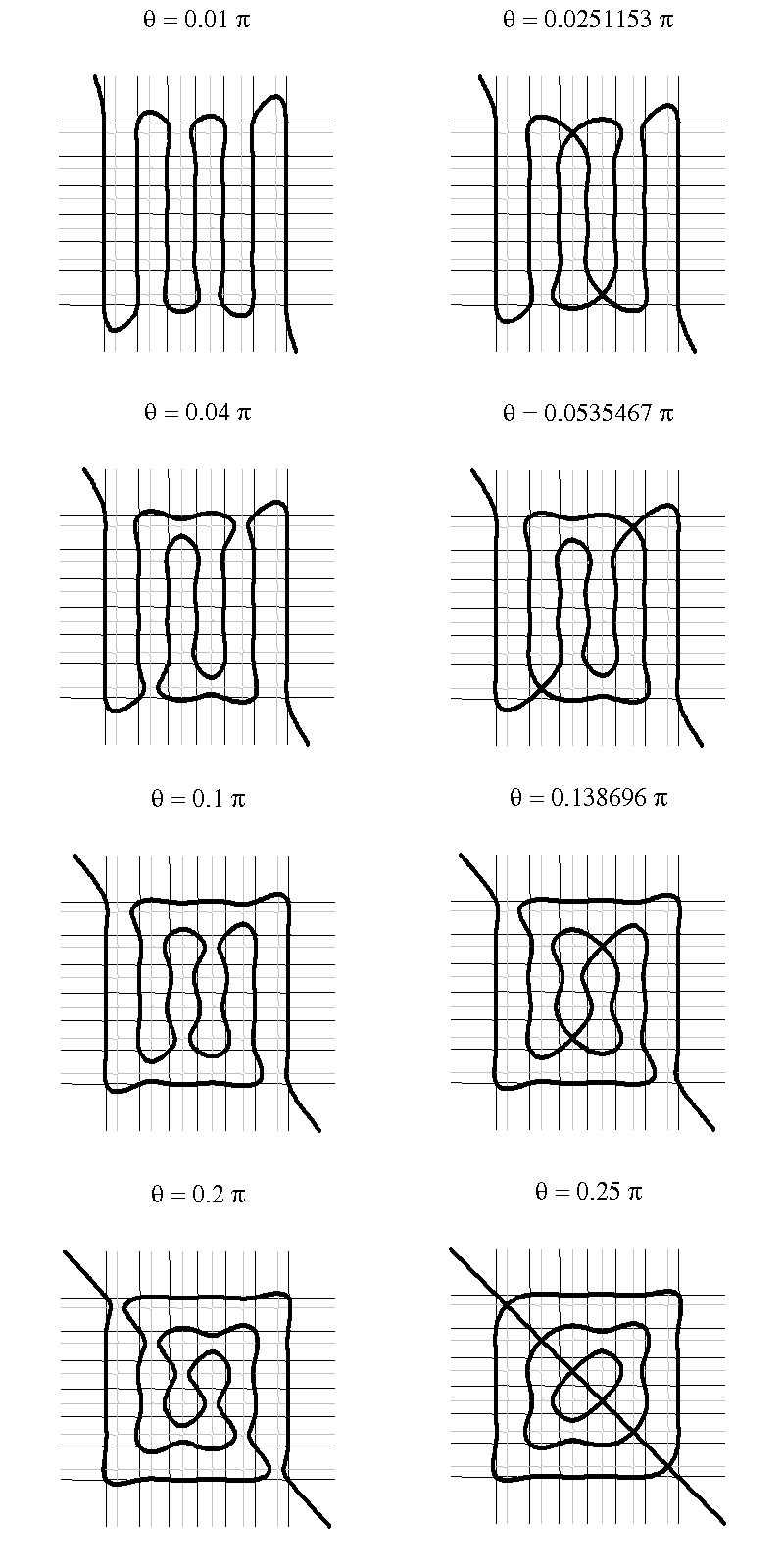}
\caption{Evolution of the nodal set $N(\Phi^{\theta}_n)$, for
$n=7$ and $\theta \in ]0,\piq]$.} \label{FSch-1}
\end{center}
\end{figure}

Figure~\ref{FSch-1}  shows the evolution of the nodal set of
$\Phi^{\theta}_n\,$,  for $n=7$, when $\theta$ varies in the interval $]0,\piq]$. The pictures in the first column correspond to regular values of $\theta$ whereas the pictures in the second column correspond to critical values of $\theta$. The form of the nodal set is stable in the open interval between two consecutive critical values of the parameter $\theta$. In the figures, the  thick curves represent the nodal sets $N(\Phi^{\theta}_7$), the thin lines correspond to the zeros of $H_7$,  and the grey lines to the zeros of $H'_7$, \ie to the zeros of $H_6$.

\medskip

We now describe the nodal set $N(\Phi^{\theta}_n)$ outside a large
enough square which contains the lattice $\cL$.  For this
purpose, we give the following two \emph{barrier lemmas} which describe the intersections of the nodal set with horizontal and vertical lines.

\begin{lemma}\label{L-ext1}
Assume that $\theta \in ]0,\piq]$.  For  $n$ odd, define
$t_{n-1,0}$ to be the unique point in $]-\infty, t_{n,1}[$ such that
$H_n(t_{n-1,0}) = - H_n(t_{n-1,1})$. Then,\vspace{-3mm}
\begin{enumerate}
    \item $\forall t \le t_{n,1}\,$, the function $y \mapsto
    \Phi^{\theta}_n(t,y)$ has exactly one zero in the interval
    $[t_{n,n},+\infty[$\,;
    \item $\forall t < t_{n-1,0}\,$, the function $y \mapsto
    \Phi^{\theta}_n(t,y)$ has exactly one zero in the interval\break
    $]-\infty,+\infty[$\,.
\end{enumerate}\vspace{-3mm}
Using the symmetry with respect to the vertical
line $\{x=0\}$, one has similar statements for $t \ge t_{n,n}$ and
for $t > - t_{n-1,0}\,$.
\end{lemma}

\textbf{Proof}.\\
Let $v(y) := \exp(\frac{t^2+y^2}{2})\, \Phi^{\theta}_n(t,y)$. In
$]t_{n,n},+\infty[\,$, $v'(y)$ is positive, and $v(t_{n,n})\le 0\,$.
The first assertion follows. The local extrema of $v$ occur at the
points $t_{n-1,j}\,$, for $1 \le j \le (n-1)\,$. The second assertion
follows from the definition of $t_{n-1,0}\,$, and from the
inequalities,
\begin{equation*}
\begin{split}
\cos\theta \, H_n(t) + \sin\theta \, H_n(t_{n-1,j}) & \le
\frac{1}{\sqrt{2}}\Big( H_n(t) + |H_n|(t_{n-1,j})\Big)\\
& < -\, \frac{1}{\sqrt{2}} \Big( H_n(t_{n-1,1}) -
|H_n|(t_{n-1,j})\Big) \le 0\,,
\end{split}
\end{equation*}
for $t < t_{n-1,0}\,$, where we have used Lemma~\ref{L-leh}.\hfill $\square$

\begin{lemma}\label{L-ext2} Let $\theta \in ]0,\piq]$. Define
$t_{n-1,n}^{\theta} \in ]t_{n,n}\,,\,\infty[$ to be the unique point
such that $\tan\theta \, H_n(t_{n-1,n}^{\theta}) = H_n(t_{n-1,1})$.
Then,\vspace{-3mm}
\begin{enumerate}
    \item $\forall t \ge t_{n,n}$, the function $x \mapsto
    \Phi^{\theta}_n(x,t)$ has exactly one zero in the interval
    $]-\infty,t_{n,1}]$\,;
    \item $\forall t > t_{n-1,n}^{\theta}$, the function $x \mapsto
    \Phi^{\theta}_n(x,t)$ has exactly one zero in the interval
    $]-\infty,\infty[$\,.
    \item For $\theta_2 > \theta_1$, we have
    $t_{n-1,n}^{\theta_2} < t_{n-1,n}^{\theta_1}\,$.
\end{enumerate}\vspace{-3mm}
Using the symmetry with respect to the horizontal
line $\{y=0\}$, one has similar statements for $t \le t_{n,1}$ and
for $t < - t_{n-1,n}^{\theta}\,$.
\end{lemma}

\textbf{Proof}. Let $h(x) := \exp(\frac{x^2+t^2}{2})\,
\Phi^{\theta}_n(x,t)$. In the interval $]-\infty,t_{n,1}]\,$, the
derivative $h'(x)$ is positive, $h(t_{n,1}) > 0\,$, and
$\lim_{x\rightarrow  -\infty} h(x) =-\infty\,$, since $n$ is odd. The
first assertion follows. The local extrema of $h$ are achieved at
the points $t_{n-1,j}$. Using Lemma~\ref{L-leh}, for $t \geq
t_{n-1,n}^{\theta}\,$, we have the inequalities,
\begin{equation*}
\begin{split}
H_n(t_{n-1,j}) + \tan\theta \, H_n(t) &\ge
\tan\theta \, H_n(t_{n-1,n}^{\theta}) - |H_n(t_{n-1,j})|\\
& = H_n(t_{n-1,1}) - |H_n(t_{n-1,j})| \ge 0\,.
\end{split}
\end{equation*}
\hfill $\square$\medskip

As a consequence of the above lemmas, we have the following
description of the nodal set far enough from $(0,0)$.

\begin{proposition}\label{P-ext}
Let $\theta \in ]0,\piq]$. In the set $\R^2 \setminus
]-t_{n-1,n}^{\theta},t_{n-1,n}^{\theta}[\times ]t_{n-1,0},
|t_{n-1,0}|[$, the nodal set $N(\Phi^{\theta}_n)$ consists of two
regular arcs. The first arc is a graph $y(x)$ over the
interval\break $]-\infty,t_{n,1}]$, starting from the point
$(t_{n,1},t_{n,n})$ and escaping to infinity with,
$$ \lim_{x\to -\infty}\frac{y(x)}{x} = - \sqrt[n]{\cot\theta}\,.$$
The second arc is the image of the first one under the symmetry with
respect to $(0,0)$ in $\R^2$.
\end{proposition}

\subsection{Local nodal patterns}\label{SS-lnp}

As in the case of the Dirichlet eigenvalues for the square, we study the possible local nodal patterns taking into account the fact that  the nodal set contains the
lattice points $\cL$, can only  visit the connected components of
the set $\{H_n(x) \, H_n(y) < 0\}$ (colored white), and consists of a simple arc at the lattice points. The following figure summarized the possible nodal patterns in the interior of the square
\cite[Figure~6.4]{BeHe},

\begin{figure}[!ht]
\begin{center}
\includegraphics[width=14cm]{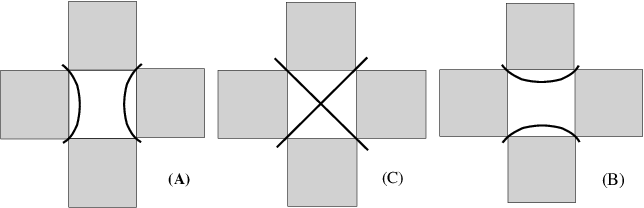}
\caption{Local nodal patterns for Dirichlet eigenfunctions of the square.} \label{FSchr-local-nodal-patterns}
\end{center}
\end{figure}

Except for nodal arcs which escape to infinity, the local nodal
patterns for the quantum harmonic oscillator are  similar (note that
in the present case, the connected components of the set $\{H_n(x)
\, H_n(y) < 0\}$ are rectangles, no longer equal squares).
{The checkerboard argument and the location of the possible critical zeros determine the possible local patterns: (A), (B) or (C). Case (C) occurs near a critical zero. Following the same ideas as in
the case of the square, in order to decide between cases (A) and
(B), we use the barrier lemmas, Lemma~\ref{L-ext1} or~\ref{L-ext2},
the vertical lines $\{x=t_{n-1,j}\}$, or the horizontal lines
$\{y=t_{n-1,j}\}$.

\section{Proof of Theorem \ref{th1}}\label{S-th1}

Note that
$$
\phi_{n,0}(x,y) - \phi_{0,n}(x,y) = - \sqrt{2}\, \Phi^{\frac{3\pi}{4}}_n(x,y) = - \sqrt{2}
\Phi^{\frac{\pi}{4}}_n(x,-y)\,.
$$

Hence, up to symmetry, it is the same to work with $\theta =\frac \pi 4$ and the anti-diagonal, or to work with $\theta = \frac{3\pi}{4}$ and the diagonal. For notational convenience, we work with $\frac{3\pi}{4}$.

\subsection{ The nodal set of $\Phi_n^{\frac{3\pi}{4}}$}

The purpose of this section is to prove the following result which is the starting point for the proof of  Theorem~\ref{th1}.

 \begin{proposition}\label{prop2} Let $\{t_{n-1,i}\,, 1 \le i \le
n-1\}$ denote the zeroes of $H_{n-1}$. For $n$ odd, the nodal set of
$\phi_{n,0} - \phi_{0,n}$ consists of the diagonal $x=y$, and of
$\frac{n-1}{2}$ disjoint simple closed curves crossing the diagonal
at the $(n-1)$ points $(t_{n-1,i}\,,t_{n-1,i})$, and the anti-diagonal
at the $(n-1)$ points $(t_{n,i}\,, - t_{n,i})$.
\end{proposition}

To prove Proposition \ref{prop2}, we first observe that it is enough
to analyze the zero set of
$$
(x,y) \mapsto \Psi_{n}(x,y)  :=H_n(x) - H_n(y)\,.
$$

\subsubsection{ Critical zeros}

The  only possible critical zeros of $\Psi_n$ are determined by
$$
H_n'(x)=0\,,\, H_n'(y)= 0\,.
$$
Hence, they consist of the
$(n-1)^2$ points $(t_{n-1,i}\,,t_{n-1,j})\,$, for  $1 \le i,j \le (n-1)\,$,
where $t_{n-1,i}$ is the $i$-th zero of the polynomial $H_{n-1}\,$.

The zero set  of $\Psi_{n}$  contains the diagonal $\{x=y\}\,$.
 Since $n$ is odd, there  are only $n$ points belonging
to the zero set on the anti-diagonal $\{x+y=0\}$.

On the diagonal, there are $(n-1)$ critical points.  We claim
that there are no critical zeros outside the diagonal.
Indeed, let $(t_{n-1,i}\,,t_{n-1,j})$ be a critical zero. Then,
$H_n(t_{n-1,i}) = H_n(t_{n-1,j})$. Using Lemma~\ref{L-leh} and the
parity properties of Hermite polynomials, we see that
$|H_n(t_{n-1,i})| = |H_n(t_{n-1,j})|$ occurs if and only if
$t_{n-1,i} = \pm t_{n-1,j}\,$. Since $n$ is odd, we can conclude that
$H_n(t_{n-1,i}) = H_n(t_{n-1,j})$ occurs if and only if
$t_{n-1,i}=t_{n-1,j}\,$.

\subsubsection{ Existence of disjoint simple closed curves in the nodal
set of $\Phi_n^{\frac{3\pi}{4}}$}

The second part in the proof of the proposition follows closely the
proof in the case of the Dirichlet Laplacian for the square (see Section 5 in \cite{BeHe}).
Essentially, the  Chebyshev polynomials are replaced by the Hermite
polynomials. Note however that the checkerboard is no more with
equal squares, and that  the square $[0,\pi]^2$ has to be
replaced in the argument by the rectangle  $[t_{n-1,0}, -t_{n-1,0}] \times
[-t_{n-1,n}^{\theta},t_{n-1,n}^{\theta}]$, for some $\theta \in ]0,\frac{3\pi}{4}[\,$, see Lemmas~\ref{L-ext1} and~\ref{L-ext2}\,.

The checkerboard argument holds, see \eqref{P24} and the properties at the beginning of Section~\ref{ss-gpns}.

The separation lemmas of our previous  paper \cite{BeHe} must
be substituted by Lemmas~\ref{L-ext1} and~\ref{L-ext2}\,, and similar
statements with the lines $\{x=t_{n-1,j}\}$ and $\{y=t_{n-1,j}\}$,
for $1 \le j \le (n-1)\,$.

One  needs to control what is going on at infinity.  As a matter of
fact, outside a specific rectangle  centered at the origin, the zero set
is the diagonal $\{x=y\}$, see Proposition~\ref{P-ext}\,.

Hence in this way (like for the square), we obtain that the nodal set of $\Psi_n$ consists of the diagonal and $\frac{n-1}{2}$  disjoint simple closed curves turning around the origin. The set $\cL$ is contained in the union of these closed curves.

\subsubsection{No other closed curve in the nodal set of
$\Phi_n^{\frac{3\pi}{4}}$}\label{nocc}

It remains to show that there are no other closed curves which do
not cross the diagonal. The ``energy'' considerations of our
previous papers \cite{BeHe, BeHe2} work here as well. Here is a simple alternative argument.


We look at the line $y=\alpha\, x $ for some $\alpha \neq 1$. The intersection of
the zero set with this line corresponds to the zeroes of the
polynomial $x\mapsto H_n(x) - H_n (\alpha\, x)$ which has at most $n$
zeroes. But in our previous construction, we get at least $n$
zeroes. So the presence of extra  curves would lead to a
contradiction for some $\alpha$. This  argument  solves the
problem at  infinity as well. \qed

\subsection{Perturbation argument}

 Figure~\ref{FSch-2} shows the desingularization of the nodal
set $N(\Phi^{\frac{3\pi}{4}}_n)$, from below and from above. The
picture is the same as in the case of the square (see Figure
\ref{sq-1}), all the critical points disappear at the same time and
in the same manner, \ie all the double crossings open up
horizontally or vertically depending whether $\theta$ is less than
or bigger than $\frac{3\pi}{4}$.

\begin{figure}[!ht]
\begin{center}
\includegraphics[width=15cm]{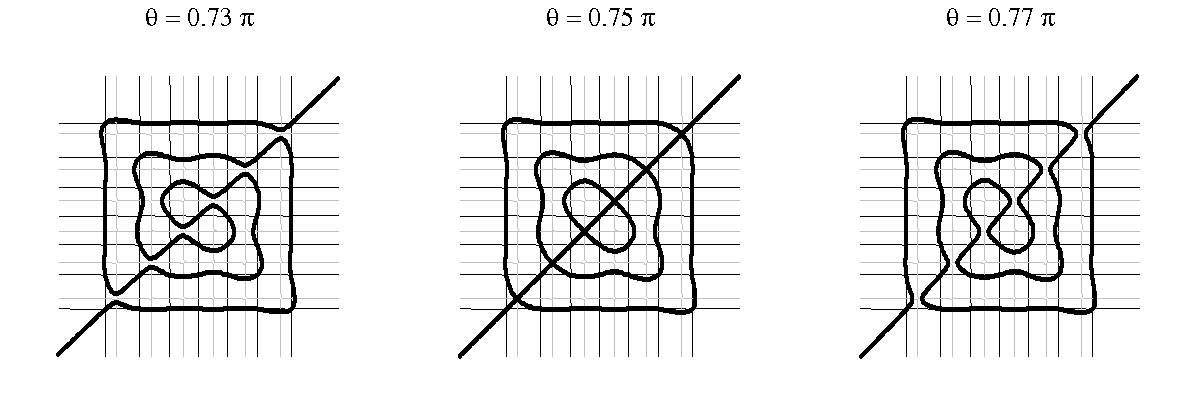}
\caption{The nodal set of $N(\Phi^{\theta}_n)$ near
$\frac{3\pi}{4}$ (here $n=7$).} \label{FSch-2}
\end{center}
\end{figure}

 As in the case of the square, in order to show that the nodal
set can be desingularized under small perturbation, we look at the
signs of the eigenfunction $\Phi_n^{\frac{3\pi}{4}}$ near the
critical zeros. We use the cases (I) and (II) which appear in Figure~\ref{FSchr-local-signs-1} below (see also
\cite[Figure~6.7]{BeHe}).

\begin{figure}[!ht]
\begin{center}
\includegraphics[width=8cm]{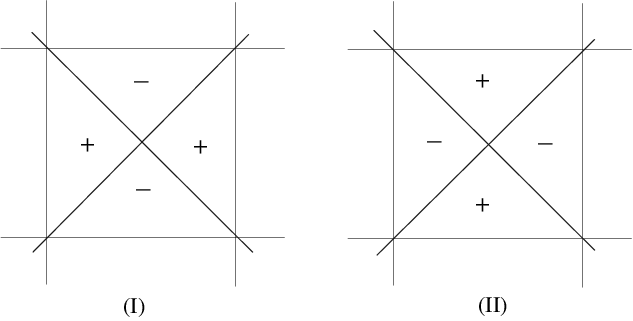}
\caption{Signs near a critical zero.} \label{FSchr-local-signs-1}
\end{center}
\end{figure}

The sign configuration for $\phi_{n,0}(x,y) - \phi_{0,n}(x,y)$ near
the critical zero $(t_{n-1,i},t_{n-1,i})$  is given by Figure~\ref{FSchr-local-signs-1}
\begin{equation*}\label{sign}
\left\lbrace
\begin{array}{l}
\text{case (I), ~if~} i \text{~is even},\\
\text{case (II),~if~} i \text{~is odd}.\\
\end{array}\right.
\end{equation*}

Looking at the intersection of the nodal set with the vertical line
$\{y=t_{n-1,i}\}$, we have that
$$(-1)^i \left( H_n(t) - H_n(t_{n-1,i})\right) \ge 0\,, \text{~for~}
t \in ]t_{n,i},t_{n,i+1}[ \,.
$$

For positive $ \epsilon$ small, we write

$$(-1)^i \left( H_n(t) - (1+\epsilon) H_n(t_{n-1,i})\right) =
(-1)^i \left( H_n(t) - H_n(t_{n-1,i})\right) + \epsilon (-1)^{i+1}\,
H_n(t_{n-1,i})\,,
$$

so that

$$(-1)^i \left( H_n(t) - (1+\epsilon) H_n(t_{n-1,i})\right) \ge 0\,,
\text{~for~} t \in ]t_{n,i},t_{n,i+1}[ \,.
$$

 A similar statement can be written for horizontal line
$\{x=t_{n-1,i}\}$ and $-\epsilon\,$, with $\epsilon > 0\,$, small
enough. These inequalities describe how the crossings  all open up  at
the same time, and in the same manner, vertically (case I) or
horizontally (case II), see Figure~\ref{FSchr-local-signs-2}, as in
the case of the square \cite[Figure~6.8]{BeHe}.

\begin{figure}[!ht]
\begin{center}
\includegraphics[width=8cm]{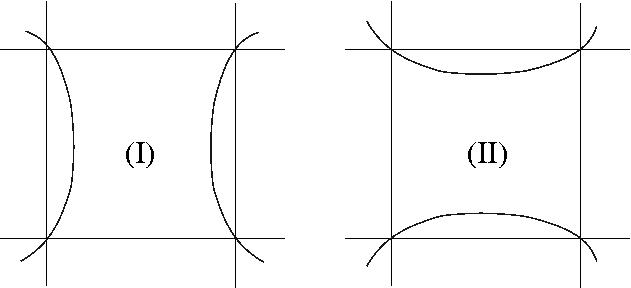}
\caption{Desingularization at a critical zero.}
\label{FSchr-local-signs-2}
\end{center}
\end{figure}

We can then conclude as in the case of the square, using the local
nodal patterns, Section~\ref{SS-lnp}.\medskip

 \textbf{Remark}. Because the local nodal patterns can only
change when $\theta$ passes through one of the values $\theta(i,j)$
defined in \eqref{P10}, the above arguments work for $\theta \in
J\setminus \{\frac{3\pi}{4}\}$, for any interval $J$ containing
$\frac{3\pi}{4}$ and no other critical value $\theta(i,j)$.

\section{Proof of Theorem \ref{th2}}\label{S-th2}

\begin{proposition}\label{P-1}
The conclusion of Theorem \ref{th2} holds with
\begin{equation}\label{P26}
\theta_c := \inf \left\{ \theta(i,j) ~|~ 1 \le i,j \le n-1 \right\}
\,,
\end{equation}
 where the critical values $\theta(i,j)$ are defined by \eqref{P10}.
\end{proposition}

\textbf{Proof}. The proof consists in the following steps. For
simplicity, we call $N$ the nodal set $N(\Phi^{\theta}_n)$.
\vspace{-3mm}
\begin{itemize}
    \item Step 1.~ By Proposition~\ref{P-ext}, the structure of the nodal
    set $N$ is known  outside a large coordinate rectangle centered at $(0,0)$
    whose sides are defined by the ad hoc numbers in
    Lemmas~\ref{L-ext1} and \ref{L-ext2}. Notice that the sides of
    the rectangle  serve as barriers for the arguments using the local
    nodal patterns as in our paper for the square.
    \item Step 2.~ For $1\le j \le n-1$, the line $\{x=t_{n-1,j}\}$
    intersects the set $N$ at exactly one point $(t_{n-1,j}, y_j)$,
    with $y_j > t_{n,n}$ when $j$ is odd, \resp with $y_j < t_{n,1}$
    when $j$ is even. The proof is given below, and is similar to the proofs of
    Lemmas~\ref{L-ext1} or \ref{L-ext2}\,.
    \item Step 3.~ Any connected component of $N$
    has at least one point in common with the set $\cL$. This
    follows from the  argument with $y = \alpha x$ or from the energy argument
    (see Subsection \ref{nocc}).
    \item Step 4.~ Follow the nodal set from the point
    $(t_{n,1},t_{n,n})$ to the point $(t_{n,n},t_{n,1})$, using the
    analysis of the local nodal patterns as in the case of the
    square.
\end{itemize}

\emph{Proof of Step~2}. For $1 \le j \le (n-1)$, define the
function $v_j$ by
$$
v_j(y) := \cos\theta \, H_n(t_{n-1,j}) + \sin\theta \, H_n(y)\,.
$$

The local extrema of $v_j$ are achieved at the points $t_{n-1,i}$,
for $1 \le i \le (n-1)$, and we have
$$
v_j(t_{n-1,i} ) = \cos\theta \, H_n(t_{n-1,j}) + \sin\theta \,
H_n(t_{n-1,i})\,,
$$
which can be rewritten, using \eqref{P10}, as
$$
v_j(t_{n-1,i}) = \frac{H_n(t_{n-1,j})}{\sin\theta(j,i)}\,
\sin\left(\theta(j,i) -\theta \right).
$$

The first term in the right-hand side has the sign of $(-1)^{j+1}$
and the second term is positive provided that $0 < \theta <
\theta_c$. Under this last assumption, we have
\begin{equation}\label{P-32}
(-1)^{j+1}\, v_j(t_{n-1,i}) > 0\,, ~\forall i\,, ~1 \le i \le (n-1)\,.
\end{equation}

The assertion follows. \hfill $\square$


\section{Eigenfunctions with ``many'' nodal domains, proof of Theorem~\ref{T-mnd}}\label{S-ovals}

This section is devoted to the proof of Theorem~\ref{T-mnd} i.e., to the constructions of eigenfunctions of $\widehat{H}$ with regular  nodal sets (no self-intersections) and ``many'' nodal domains. We work in polar coordinates. An orthogonal basis of $\cE_{\ell}$ is given by the functions $\Omega^{\pm}_{\ell,n}$,
\begin{equation}\label{ov6}
\Omega^{\pm}_{\ell,n}(r,\varphi ) =  \exp(-\frac{r^2}{2}) \,
r^{\ell-2n}\, L_{n}^{(\ell-2n)}(r^2)\, \exp \left( \pm i
(\ell-2n)\varphi  \right) \,,
\end{equation}
with $0 \le n \le \left[ \frac{\ell}{2} \right] $, see \cite[Section~2.1]{LeyD}. In this formula $L_n^{(\alpha)}$ is the generalized Laguerre polynomial of degree $n$ and parameter $\alpha$, see \cite[Chapter~5]{Sz}. Recall that the Laguerre polynomial $L_n$ is the polynomial $L_n^{(0)}$.

\begin{assumption}\label{a-ovals}
From now on, we assume that $\ell = 4k\,$, with $k$ even.
\end{assumption}%

Since $\ell$ is even, we have a rotation invariant eigenfunction $\exp(-\frac{r^2}{2}) \, L_{2k}(r^2)$ which has $(2k +1)$ nodal
domains. We also look at the eigenfunctions $\omega_{\ell,n}$,
\begin{equation}\label{ov8}
\omega_{\ell,n}(r,\varphi ) =  \exp(-\frac{r^2}{2})  \,
r^{\ell-2n}\, L_{n}^{(\ell-2n)}(r^2)\, \sin\left((\ell-2n)\varphi
\right)\,,
\end{equation}
with $0 \le n < \left[\frac{\ell}{2}\right]$. \\
The number of nodal domains of these eigenfunctions is $\mu(\omega_{\ell,n}) = 2(n+1)(\ell-2n)$, because the Laguerre polynomial of degree $n$ has $n$ simple positive roots. When $\ell = 4k$,  the largest of these numbers is
\begin{equation}\label{ov10}
\mu_\ell : = 4k(k+1)\,,
\end{equation}
and this is achieved for $n=k\,$.

When $k$ tends to infinity, we have $ \mu_{\ell} \sim \frac{\ell^2}{4}$, the same order of magnitude as Leydold's upper bound $ \mu_\ell^L  \sim \frac{\ell^2}{2}$.

We want now to construct eigenfunctions $u_k \in \cE_{\ell}\,, \ell = 4k$ with \emph{regular} nodal sets, and ``many'' nodal domains (or equivalently, ``many'' nodal connected components), more precisely with $\mu(u_k) \sim \frac{\ell^2}{8}$.

The construction consists of the following steps. \vspace{-3mm}
\begin{enumerate}
    \item Choose $A \in \cE_{\ell}$ such that $\mu(A) = \mu_\ell \,$.
    \item Choose $B \in \cE_{\ell}$ such that for $a$ small enough,
    the perturbed eigenfunction \break $F_{a}:=A+a\,B$ has no critical zero
    except the origin, and a nodal set with many components.
    Fix such an $a$\,.
    \item Choose $C \in \cE_{\ell}$ such that for $b$ small enough
    (and $a$ fixed), $G_{a,b}:=A+a\,B+b\,C$ has no critical zero.
\end{enumerate}

From now on, we fix some $\epsilon$, $0 < \epsilon < 1\,$. We assume
that $a$ is positive (to be chosen small enough later on), and
that $b$ is non zero (to be chosen small enough, either positive or
negative later on). \bigskip

In the remaining part of this section, we skip the exponential
factor in the eigenfunctions since it is irrelevant to study the
nodal sets.\bigskip

Under Assumption~\ref{a-ovals}, define
\begin{equation}\label{ov12}
A(r,\varphi ) := r^{2k}\, L_k^{(2k)}(r^2) \, \sin(2k\varphi )\,,\,
B(r,\varphi ) := r^{4k} \, \sin(4k\varphi  - \epsilon \pi)\,,\,
C(r,\varphi ) := L_{2k}(r^2)\,.
\end{equation}

We consider the deformations $F_{a} = A+a\,B$ and $G_{a,b} = A+a\,B+b\,C$.
Both functions are invariant under the rotation of angle
$\frac{\pi}{k}$, so that we can restrict to $\varphi \in
[0,\frac{\pi}{k}]\,$.

For later purposes, we introduce the angles $\varphi_j =
\frac{j\pi}{k}$, for $0 \le j \le 4k-1$, and $\psi_m =
\frac{(m+\epsilon)\pi}{4k}$, for $0 \le m \le 8k-1$\,. We denote by
$t_i, 1 \le i \le k$ the zeros of $L_k^{(2k)}$, listed in increasing
order. They are simple and positive, so that the numbers $r_i =
\sqrt{t_i}$ are well defined. For notational convenience, we denote by $\dot{L}_{k}^{(2k)}$ the
derivative of the polynomial $L_k^{(2k)}$. This polynomial has $(k-1)$ simple
zeros, which we denote by $t_i', 1 \le i \le k-1$, with $t_i < t_i' < t_{i+1}\,$. We
define $r_i' := \sqrt{t_i'}\,$.

\subsection{Critical zeros}\label{SS-cz}

Clearly, the origin is a critical zero of the eigenfunction $F_{a} =
A+a\,B\,$, while \break  $G_{a,b} = A+a\,B+b\,C$ does not vanish at the origin.

Away from the origin, the critical zeros of $F_{a}$ are given by the
system
\begin{equation}\label{ov14}
\begin{array}{ll}
F_{a}(r,\varphi ) &= 0\,,\\
\partial_{r}F_{a}(r,\varphi ) &= 0\,,\\
\partial_{\varphi }F_{a}(r,\varphi ) &= 0\,.\\
\end{array}
\end{equation}

The first and second conditions imply that a critical zero
$(r,\varphi )$ satisfies

\begin{equation}\label{ov14a}
\sin(2k\varphi ) \, \sin(4k\varphi -\epsilon \pi)\, \left( k
L_k^{(2k)}(r^2) - r^2 \dot{L}_{k}^{(2k)}(r^2)\right) = 0\,,
\end{equation}
where $\dot{L}$ is the derivative of the polynomial $L\,$.

The first and third conditions imply that a critical zero
$(r,\varphi )$ satisfies

\begin{equation}\label{ov14b}
2 \sin(2k\varphi ) \, \cos(4k\varphi  -\epsilon \pi) -
\cos(2k\varphi ) \, \sin(4k\varphi  -\epsilon \pi) = 0\,.
\end{equation}

It is easy to deduce from \eqref{ov14} that when $(r,\varphi)$ is a
critical zero,\break  $\sin(2k\varphi ) \, \sin(4k\varphi -\epsilon \pi)
\neq 0\,$. It follows that, away from the origin, a critical zero
$(r,\varphi )$ of $F_{a}$ satisfies the system
\begin{equation}\label{ov16}
\begin{array}{ll}
k L_k^{(2k)}(r^2) - r^2 \dot{L}_{k}^{(2k)}(r^2) = 0\,,\\[5pt]
2 \sin(2k\varphi ) \, \cos(4k\varphi  -\epsilon \pi) -
\cos(2k\varphi ) \, \sin(4k\varphi  -\epsilon \pi) = 0\,.
\end{array}
\end{equation}

The first equation has precisely $(k-1)$ positive simple zeros
$r_{c,i}$, one in each interval $]r_i,r'_i[\,$, for $1 \le i \le
(k-1)$. An easy analysis of the second shows that it has $4k$ simple
zeros $\varphi_{c,j}$, one in each interval
$]\varphi_j,\psi_{2j+1}[\,$, for $0 \le j \le 4k-1\,$.

\begin{property}\label{P-ov2}
The only possible critical zeros of the function $F_{a}$, away from
the origin, are the points $(r_{c,i},\varphi_{c,j})\,$, for $1 \le i
\le k-1$ and $0 \le j \le 4k-1\,$, with corresponding finitely many
values of $a$ given by \eqref{ov14}. In particular, there exists some
$a_0 > 0$ such that for $0 < a < a_0$, the eigenfunction $F_{a}$ has
no critical zero away from the origin.
\end{property}

The function $G_{a,b}$ does not vanish at the origin (provided that
$b \neq 0$). Its critical zeros are given by the system

\begin{equation}\label{ov20}
\begin{array}{ll}
G(r,\varphi) & = 0\,,\\
\partial_r G(r,\varphi) & = 0\,,\\
\partial_{\varphi} G(r,\varphi) & = 0\,.\\
\end{array}
\end{equation}

We look at the situation for $r$ large. Write

\begin{equation}\label{ov22}
\begin{array}{ll}
L_k^{(2k)}(t) & = \frac{(-1)^k}{k!}\, t^k + P_k(t)\,,\\[5pt]
L_{2k}(t) & = \frac{1}{(2k)!} \, t^{2k} + Q_k(t)\,,\\
\end{array}
\end{equation}
where $P_k$ and $Q_k$ are polynomials with degree $(k-1)$ and
$(2k-1)$ respectively.\medskip

The first and second equations in \eqref{ov20} are equivalent to the
first and second equations of the system
\begin{equation}\label{ov24}
\begin{array}{ll}
0 & = \frac{(-1)^k}{k!} \, \sin(2k\varphi) + a \, \sin(4k\varphi
-\epsilon \pi) + \frac{b}{(2k)!} + O(\frac{1}{r^2})\,,\\
0 & = \frac{(-1)^k}{k!}\, \cos(2k\varphi) + 2 a \cos(4k\varphi
-\epsilon \pi) + O(\frac{1}{r^2})\,,\\
\end{array}
\end{equation}
where the $O(\frac{1}{r^2})$ are uniform in $\varphi$ and $a, b$
(provided they are initially bounded).

\begin{property}\label{P-ov4}
There exist positive numbers $a_1 \le a_0$, $b_1, R_1$, such that
for $0 < a < a_1$, $0 < |b| < b_1$, and $r > R_1$, the function
$G_{a,b}(r,\varphi)$ has no critical zero. It follows that for fixed
$0 < a < a_1$, and $b$ small enough (depending on $a$), the function
$G_{a,b}$ has no critical zero in $\R^2$.
\end{property}

\textbf{Proof}. Let $\alpha := \frac{(-1)^k}{k!} \, \sin(2k\varphi)
+ a \, \sin(4k\varphi -\epsilon \pi)$, $\beta := \frac{b}{(2k)!}$
and $$\gamma := \frac{(-1)^k}{k!}\, \cos(2k\varphi) + 2 a
\cos(4k\varphi -\epsilon \pi)\,.$$
 Compute $(\alpha + \beta)^2 +
\gamma^2$. For $0 < a < \frac{1}{2 k!}\,$, one has
$$
(\alpha + \beta)^2 + \gamma^2 \le \frac{1}{(2 k!)^2} - \frac{4a}{k!}
- \frac{4|b|}{k! (2k)!}\,.
$$
The first assertion follows.  The second assertion follows from the
first one and from Property~\ref{P-ov2}. \hfill $\square$

\subsection{The checkerboard}\label{SS-cb}

Since $a$ in positive, the nodal set of $F_{a}$ satisfies
\begin{equation}\label{cb2}
\cL \subset N(F_{a}) \subset \cL \cup \left\lbrace A B < 0
\right\rbrace ,
\end{equation}
where $\cL$ is the finite set $N(A)\cap N(B)$, more precisely,
\begin{equation}\label{cb4}
\cL = \left\lbrace (r_i,\psi_m) ~|~ 1 \le i \le k\,, ~ 0 \le m \le
8k-1 \right\rbrace .
\end{equation}

Let $p_{i,m}$ denote the point with polar coordinates
$(r_i,\psi_m)$. It is easy to check that the points $p_{i,m}$ are
regular points of the nodal set $N(F_{a})$. More precisely the nodal
set $N(F_{a})$ at these points is a regular arc transversal to the
lines $\{\varphi = \psi_m\}$ and $\{r=r_i\}$. Note also that the
nodal set $N(F_{a})$ can only cross the nodal sets $N(A)$ or $N(B)$
at the points in $\cL\,$. \medskip

The connected components of the set $\left\lbrace A B \neq 0
\right\rbrace$ form a ``polar checkerboard'' whose white boxes are
the connected components in which $AB < 0\,$. The global aspect of the
checkerboard depends on the parity of $k\,$. Recall that our assumption is that $\ell = 4k\,$, with $k$ even. Figure~\ref{F-cbe} displays a partial view of the checkerboard, using the invariance under the
rotation of angle $\frac{\pi}{k}$. The thin lines labelled ``R'' correspond to the angles $\psi_m$, with $m=0, 1, 2, 3\,$. The thick lines to the angles $\varphi_j$, with $j=0,1, 2\,$. The thick arcs of circle correspond to the values $r_i\,$, with $i=1, 2, 3$ and then $i=k-1, k\,$. The light grey part
represents the zone $r_i$ with $i=4, \ldots, k-2$.  The intersection points of the thin lines ``R'' with the thick arcs are the point in $\cL\,$, in the sector $0 \le \varphi \le \varphi_2\,$. The outer arc of circle (in grey) represents the horizon.

\begin{figure}
\begin{center}
\includegraphics[width=7cm]{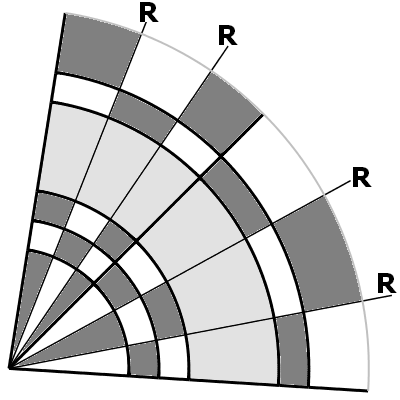}
\caption{$\ell = 4k$, $k$ even.}\label{F-cbe}
\end{center}
\end{figure}

\subsection{Behavior at infinity}\label{SS-bi}

We now look at the behavior at infinity of the functions $F_{a}$
and $G_{a,b}\,$. We restrict our attention to the sector $\{0 \le
\varphi \le \varphi_2\}$.

Recall that $k$ is even.

For $r > r_k$, the nodal set $N(F_{a})$ can only visit the white
sectors $\cS_0 := \{\varphi_0 < \varphi < \psi_0\}$, $\cS_1 :=
\{\psi_1 <\varphi < \varphi_1\}$, and $\cS_2 :=\{ \psi_2 < \varphi <
\psi_3 \}$, issuing respectively from the points $p_{k,0}, p_{k,1}$
or $p_{k,2}, p_{k,3}\,$.

As above, we can write
\begin{equation}\label{bi40}
F_{a}(r,\varphi) = r^{4k} g(\varphi) + \sin(2k\varphi) P_k(r^2)\,,
\end{equation}
with $$g(\varphi) = \frac{1}{k!}\sin(2k\varphi) + a \sin(4k\varphi
-\epsilon \pi)\,,
$$
where we have used the fact that $k$ is even.\medskip

\noib \textbf{Analysis in} $\cS_0$.~We have $0 < \varphi <
\frac{\epsilon \pi}{4k}\,$. Note that $g(0) \, g(\frac{\epsilon \pi}{4k})
< 0\,$. On the other-hand, $g'(\varphi)$ satisfies
$$
g'(\varphi) \ge 2k \left\lbrace \frac{1}{k!} \cos(\frac{\epsilon
\pi}{2}) - 2a\right\rbrace .
$$

It follows that provided that $0 < a < \frac{1}{2\,k!}
\cos(\frac{\epsilon \pi}{2})$, the function $g$ has exactly one zero
$\theta_0$ in the interval $]0,\frac{\epsilon \pi}{4k}[\,$.

It follows that for $r$ big enough, the equation $F_{a}(r,\varphi) =
0$ has exactly one zero $\varphi(r)$ in the interval $]0,\frac{\epsilon
\pi}{4k}[$, and this zero tends to $\theta_0$ when $r$ tends to
infinity. Looking at \eqref{bi40} again, we see that $\varphi(r) =
\theta_0 + O(\frac{1}{r^2})$. It follows that the nodal set in the
sector $\cS_0$ is a line issuing from $p_{k,0}$ and tending to
infinity with the asymptote $\varphi = \theta_0\,$.\bigskip

\noib \textbf{Analysis in} $\cS_1$.~The analysis is similar to the
analysis in $\cS_0$.\bigskip

\noib \textbf{Analysis in} $\cS_2$.~In this case, we have that
$\frac{(2+\epsilon)\pi}{4k} < \varphi < \frac{(3+\epsilon)\pi}{4k}\,$.
It follows that $- \sin(2k\varphi) \ge \min\{ \sin(\frac{\epsilon
\pi}{2}), \cos(\frac{\epsilon \pi}{2})\} > 0$. If $0 < a < \frac{1}{k!} \min\{ \sin(\frac{\epsilon \pi}{2}),
\cos(\frac{\epsilon \pi}{2})\}$, then $F_{a}(r,\varphi)$ tends to
negative infinity when $r$ tends to infinity, uniformly in $\varphi
\in ]\psi_2,\psi_3[\,$. It follows that the nodal set of $N(F_{a})$ is
bounded in the sector $\cS_2\,$.

\subsection{The nodal set $N(F_{a})$ and $N(G_{a,b})$}\label{SS-nodalFa}

\begin{proposition}\label{P-Fa-2}
For $\ell = 4k$, $k$ even, and $a$ positive small enough, the nodal
set of $F_{a}$ consists of three sets of ``ovals''\vspace{-3mm}
\begin{enumerate}
\item a cluster of $2k$ closed (singular) curves, with a common singular
point at the origin,
\item $2k$ curves going to infinity, tangentially to
lines $\varphi  = \vartheta_{a,j}$ (in the case of the sphere, they
would correspond to a cluster of closed curves at the south pole),
\item $2k(k-1)$ disjoint simple closed curves (which correspond to
the white cases at finite distance of Stern's checkerboard for $A$
and $B$).
\end{enumerate}
\end{proposition}

\textbf{Proof}.

Since $B$ vanishes at higher order than $A$ at the origin, the
behavior of the nodal set of $F_{a}$ is well determined at the
origin. More precisely, the nodal set of $F_{a}$ at the origin
consists of $4k$ semi-arcs, issuing from the origin tangentially to
the lines $\varphi  = j\pi/2k$, for $0 \le j \le 4k-1\,$. At infinity,
the behavior of the nodal set of $F_{a}$ is determined for $a$
small enough in Subsection~\ref{SS-bi}. An analysis \`{a} la Stern, then
shows that for $a$ small enough there is a cluster of  ovals  in the
intermediate region $\{r_2 < r < r_{k-1}\}$, when $k \ge 4$. \hfill $\square$

Fixing so $a$ small enough so that the preceding proposition holds, in order to obtain a regular nodal set, it suffices to perturb $F_a$ into
$G_{a,b}\,$, with $b$ small enough, choosing its sign so that the nodal set $F_{a}$ is desingularized at the origin, creating $2k$ ovals.

Figure~\ref{F-ov4K-K2} displays the cases $\ell = 8$ (\ie $k=2$).

\begin{figure}[!ht]
\begin{center}
\includegraphics[width=15cm]{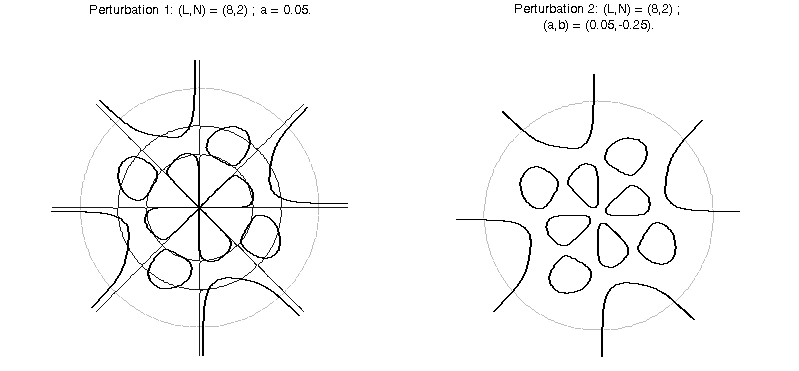}
\caption{Ovals for $\ell = 8\,$.} \label{F-ov4K-K2}
\end{center}
\end{figure}

Finally, we have constructed an eigenfunction $G_{a,b}$ with $2k(k+1)$ nodal component so that $\mu(G_{a,b}) \sim 2 k^2 = \frac{\ell^2}{8}\,$.


\section{On bounds for the length of the nodal set}\label{S-lbnz}

In Subsection~\ref{ss-lbbg}, we obtain Theorem~\ref{lbnz-P1} as a corollary of a more general result, Proposition~\ref{lbnz-P2}, which sheds some light on the exponent $\frac 32$.  The proof is typically $2$-dimensional, \`{a} la Br\"{u}ning-Gromes \cite{BrGr,Br78}. We consider more general potentials in Subsection~\ref{SS-lj}. In Subsection~\ref{SS-lj}, we extend the methods of Long Jin \cite{LJ} to some Schr\"{o}dinger operators. We obtain both lower and upper bounds on the length of the nodal sets in the classically permitted region, Proposition~\ref{lje}.

\subsection{Lower bounds, proof \`{a} la Br\"{u}ning-Gromes} \label{ss-lbbg}


Consider the eigenvalue problem on $L^2(\R^2)$
\begin{equation}\label{lbnz-4}
H_V := - \Delta + V(x)\,, ~~H_V\, u = \lambda\, u\,,
\end{equation}
for some suitable non-negative potential $V$ such that the operator has discrete spectrum (see  \cite[Chapter~8]{Lau}).
More precisely, we assume:

\begin{assumption}\label{assump}
The potential  $V$ is positive, continuous  and tends to infinity at infinity.
 \end{assumption}

Introduce the sets
\begin{equation}\label{lbnz-6}
B_V(\lambda) := \left\lbrace  x \in \R^2 ~|~ V(x) < \lambda \right\rbrace\,,
\end{equation}
and, for $r > 0\,$,
\begin{equation}\label{lbnz-8}
B_V^{(-r)}(\lambda) := \left\lbrace x\in \R^2 ~|~ B(x,r) \subset B_V(\lambda)  \right\rbrace,
\end{equation}
where $B(x,r)$ is the open ball with center $x$ and radius $r$.

\begin{proposition}\label{lbnz-P2}
Fix $\delta \in ]0,1[$ and $\rho \in ]0,1]$. Under Assumption~\ref{assump}, for $\lambda$ large enough, and for any nonzero eigenfunction $u$ of $H_V$, $H_V u = \lambda u$, the length of $N(u) \cap B_V(\delta\lambda)$ is not less than
\begin{equation}\label{lbnz-16}
\frac{2(1-\delta)}{9\pi^2j_{0,1}}\sqrt{\lambda} \, A\left( B_V^{(-2\rho)}(\delta\lambda)\right).
\end{equation}
\end{proposition}%

\textbf{Proof} of Proposition~\ref{lbnz-P2}.

\begin{lemma}\label{lbnz-L1}
Choose some radius $0< \rho \le 1$, and let
\begin{equation}\label{lbnz-10}
\rho_{\delta} := \frac{j_{0,1}}{\sqrt{1-\delta}}\,.
\end{equation}
Then, for $ \lambda > \left(\frac{\rho_{\delta}}{\rho}\right)^2$, and for any $x \in B_V^{(-\rho)}(\delta\lambda)$, the ball $B(x, \frac{\rho_{\delta}}{\sqrt{\lambda}})$ intersects the nodal set $N(u)$ of the function $u$.
\end{lemma}%

\emph{Proof of Lemma~\ref{lbnz-L1}.}  Let $r:= \frac{\rho_{\delta}}{\sqrt{\lambda}}\,$.  If the ball $B(x,r)$ did not intersect $N(u)$, then it would be contained in a nodal domain $D$ of the eigenfunction $u$. Denoting by $\sigma_1(\Omega)$ the least Dirichlet eigenvalue of the operator $H_V$ in the domain $\Omega$, by monotonicity, we could write
$$
\lambda = \sigma_1(D) \le \sigma_1\left(B(x,r)\right) .
$$
Since $x \in B_V^{(-\rho)}(\delta\lambda)$ and $ \lambda > \left(\frac{\rho_{\delta}}{\rho}\right)^2\,$, the ball $B(x,r)$ is contained in $B_V(\delta\lambda)$, and we can bound $V$ from above by $\delta\lambda$ in this ball. It follows that $\sigma_1\left( B(x,r) \right) < \frac{j_{0,1}^2}{r^2} + \delta\lambda$. This leads to a contradiction with the definition of $\rho_{\delta}$. \hfill \qed

Consider the set $\cF$ of finite subsets $\{x_1, \ldots, x_n\}$ of $\R^2$ with the following properties,
\begin{equation}\label{lbnz-12}
\left\{
\begin{array}{l}
x_i \in N(u) \cap B_V^{(-\rho)}(\delta\lambda), ~ 1\le i \le n\,,\\[5pt]
B(x_i,\frac{\rho_{\delta}}{\sqrt{\lambda}})\,, ~ 1 \le i \le n\,, \text{~pairwise disjoint.}
\end{array}
\right.
\end{equation}

For $\lambda$ large enough, the set $\cF$ is not empty, and can be ordered by inclusion. It admits a maximal element $\{x_1, \ldots, x_N\}\,$, where $N$ depends on $\delta, \rho, \lambda$ and $u$.

\begin{lemma}\label{lbnz-L2}
The balls $B(x_i, \frac{3\rho_{\delta}}{\sqrt{\lambda}})\,, 1 \le i \le N\,$, cover the set $B_V^{(-2\rho)}(\delta\lambda)\,$.
\end{lemma}%

\emph{Proof of Lemma~\ref{lbnz-L2}.} Assume the claim in not true, \ie that there exists some $y \in B_V^{(-2\rho)}(\delta\lambda)$ such that $|y-x_i|> \frac{3\rho_{\delta}}{\sqrt{\lambda}}$ for all $i \in \{1, \ldots, N\}$. Since $y \in B_V^{(-\rho)}(\delta\lambda)\,$, by Lemma~\ref{lbnz-L1}, there exists some $x \in N(u) \cap B(y, \frac{\rho_{\delta}}{\sqrt{\lambda}})$, and we have $x \in B_V^{(-\rho)}(\delta\lambda)$. Furthermore, for all $i \in \{1, \ldots, N\}$, we have $|x-x_i| \ge \frac{2\rho_{\delta}}{\sqrt{\lambda}}\,$. The set $\{x, x_1, x_2, \ldots, x_N\}$ would belong to $\cF$, contradicting the maximality of $\{x_1, x_2, \ldots, x_N\}$. \qed

Lemma~\ref{lbnz-L2} gives a lower bound on the number $N$,
\begin{equation}\label{lbnz-14}
N \ge \frac{\lambda}{9\pi^2 \rho_{\delta}^2} A\left( B_V^{(-2\rho)}(\delta\lambda) \right)\,,
\end{equation}
where $A(\Omega)$ denotes the area of the set $\Omega\,$.

\begin{lemma}\label{lbnz-L3}
For any $\alpha < j_{0,1}\,$, the ball $B(x,\frac{\alpha}{\sqrt{\lambda}})$ does not contain any closed connected component of the nodal set $N(u)$.
\end{lemma}%

\emph{Proof of  Lemma~\ref{lbnz-L3}.} Indeed, any closed connected component of $N(u)$ contained in $B(x,\frac{\alpha}{\sqrt{\lambda}})$ would bound some nodal domain $D$ of $u$, contained in $B(x,\frac{\alpha}{\sqrt{\lambda}})$, and we would have
$$
\lambda = \sigma_1(D) \ge \sigma_1\left( B(x,\frac{\alpha}{\sqrt{\lambda}} \right) \ge \frac{j_{0,1}^2}{\alpha^2}\lambda\,,
$$
contradicting the assumption on $\alpha\,$.\quad \qed

Take the maximal set $\{x_1, \ldots, x_N\} \subset N(u) \cap B_V^{(-\rho)}(\delta\lambda)$ constructed above. The balls $B(x_i,\frac{\rho_{\delta}}{\sqrt{\lambda}})$ are pairwise disjoint, and so are the balls $B(x_i,\frac{\alpha}{\sqrt{\lambda}})$ for any $0 <\alpha < j_{0,1}$. There are at least two nodal arcs issuing from a point $x_i$, and they must exit $B(x_i, \frac{\alpha}{\sqrt{\lambda}})$, otherwise we could find a closed connected component of $N(u)$ inside this ball, contradicting Lemma~\ref{lbnz-L3}. The length of $N(u) \cap B(x_i,\frac{\alpha}{\sqrt{\lambda}})$ is at least $\frac{2\alpha}{\sqrt{\lambda}}$. Finally, the length of $N(u) \cap B_V(\delta\lambda)$ is at least $N \frac{2\alpha}{\sqrt{\lambda}}$ which is bigger than
$$\frac{2\alpha}{9\pi^2\rho_{\delta}^2} \sqrt{\lambda} A\left( B_V^{(-2\rho)}(\delta\lambda)\right)\,.
$$
Since this is true for any $\alpha < j_{0,1}$. \hfill \qed \medskip

\textbf{Proof of Theorem~\ref{lbnz-P1}.}

We apply the preceding proposition with $V(x) = |x|^{2k}$ and $\rho=1$. Then, $B_V(\lambda) = B( \lambda^{\frac{1}{2k}})$ and $B_V^{(-r)}(\lambda) = B( \lambda^{\frac{1}{2k}}-r)$. In this case, the length of the nodal set is bounded by some constant times $\lambda^{\frac{1}{2} + \frac{1}{k}} \approx  \lambda^{\frac{1}{2}}\, A\left( B_V(\delta\lambda)\right) $. When $k=1$, we obtain Proposition~\ref{lbnz-P1}. \hfill \qed \medskip

\textbf{Remark}. The above proof sheds some light on the exponent $\frac{3}{2}$ in Proposition~\ref{lbnz-P1}.

\subsection{ More general potentials}\label{S-lbnsV}

We  reinterpret Proposition~\ref{lbnz-P2}  for  more general potentials $V(x)$, under natural assumptions which appear in the determination of the Weyl's asymptotics of $H_V$ {  (see \cite{Rob}, \cite{Ho}).}\\
After renormalization, we assume:
\begin{assumption}
$V$ is of class $C^1$,
$
V\geq 1\,,$
and
there exist some positive constants $\rho_0$ and $C_1$ such that for all $x \in \R^2$,
\begin{equation}\label{v3}
|\nabla V(x)| \leq C_1  V(x)^{1-\rho_0}\,.
\end{equation}
\end{assumption}

Note that under this assumption there exist positive constants $r_0$ and $C_0$ such that
\begin{equation}\label{v4}
x,y \mbox{ satisfy } |x-y| \leq r_0 \Rightarrow V(x) \leq C_0 \, V(y)\,.
\end{equation}
The proof is easy. We first write
$$
V(x) \leq V (y) + |x-y| \sup_{z\in [x,y]} |\nabla V (z) |\,.
$$
Applying  \eqref{v3} (here we only use $\rho_0\geq 0$), we get
$$
V(x) \leq V(y) + C_1 |x-y| \sup_{z\in [x,y]} V (z) \,.
$$
We now take $x\in B(y,r)$ for some $r>0$  and get
$$
\sup_{x\in B(y,r)} V(x) \leq V(y) + C_1\,  r \, \sup_{x\in B(y,r)}  V(x)\,,
$$
which we can rewrite, if $C_1 r < 1$, in the form
$$
V(y) \leq \sup_{x\in B(y,r)} V(x) \leq V(y) (1- C_1 r)^{-1}\,.
$$
This is more precise  than \eqref{v4} because we get $C_0(r_0) = (1- C_1 r_0)^{-1}$, which tends to $1$ as $r_0 \rightarrow 0\,$.

We assume
\begin{assumption}
For any $\delta \in ]0,1[$, there exists some positive  constants $A_{\delta}$ and $\lambda_\delta$ such that
\begin{equation}\label{v5}
1 < A(B_V(\lambda))/A(B_V(\delta \lambda)) \leq A_{\delta}\,,\, \forall \lambda \geq \lambda_\delta\,.
\end{equation}
\end{assumption}

\begin{proposition}\label{prop7.9}
Fix $\delta \in (0,1)$, and assume that $V$ satisfies the previous assumptions.  Then, there exists a positive constant $C_\delta$ (depending only on the constants appearing in the assumptions on $V$) and $\lambda_\delta$ such that for any eigenpair $(u,\lambda)$ of $H_V$ with $\lambda \geq \lambda_\delta $,  the length of $ N(u) \cap B_V (\delta \lambda) $ is larger than $C_\delta \lambda^\frac 12  A(B_V(\lambda))$.
\end{proposition}

{\bf Proof}.

Using \eqref{v5}, it is enough to prove the existence of $r_1$ such that, for $0<r<r_1$, there exists $C_2(r)$ and $M(r)$ s.t.
$$
 B_V (\mu - C_2 \mu^{1-\rho_0}) \subset B_V^{-r} (\mu)\,,\, \forall \mu > M(r)\,.
$$
But, if $ x\in  B_V (\mu - C_2\mu^{1-\rho_0})$, and $y \in B (x,r)$, we have
$$
V(y) \leq V (x) + C_1  C_0(r) ^{1-\rho_0} r V(x)^{1-\rho_0} \leq \mu - C_2 \mu^{1-\rho_0} +  C_1  C_0(r) ^{1-\rho_0} r \mu^{1-\rho_0}\,.
$$
Taking  $C_2 (r)=C_1 C_0(r)^{1-\rho_0} r$  and $M(r) \geq (C_2(r)+1)^{\frac {1}{\rho_0}} $ gives the result. \qed

\textbf{Remarks.}\vspace{-3mm}\begin{enumerate}
\item The method of proof of Proposition~\ref{lbnz-P2}, which is reminiscent of the proof by Br\"{u}ning  \cite{Br78} (see also  \cite{BrGr})  is typically $2$-dimensional.
\item  The same method could be applied to a Schr\"{o}dinger operator on a complete noncompact Riemannian surface, provided one has some control on the geometry, the  first eigenvalue of small balls, etc..
\item If we assume that there exist positive constants $m_0 \leq m_1$ and $C_3$ such that for any $x \in \R^2$,
\begin{equation}\label{v2}
\frac{1}{C_3} <x>^{m_0} \leq V(x) \leq C_3 <x>^{m_1}\,,
\end{equation}
where $<x> := \sqrt{ 1+|x|^2} $, then $A(B_V(\lambda))$ has a controlled growth at $\infty$.
\item If $m_0=m_1$ in \eqref{v2}, then \eqref{v5} is satisfied. The control of $A_\delta$ as $\delta \rightarrow +1$ can be obtained under additional assumptions.
\end{enumerate}

\subsection{Upper and lower bounds on the length of the nodal set: the semi-classical approach of Long Jin}\label{SS-lj}

In \cite{LJ}, Long Jin analyzes the same question in the semi-classical context for a Schr\"odinger operator
$$H_{W,h}:=-h^2  \Delta_g + W(x)\,,
$$
where $\Delta_g$ is the Laplace-Beltrami operator on the compact connected analytic Riemannian surface $(M,g)$, with $W$ analytic.   In this context, he shows that if  $(u_h,\lambda_h)$  is an $h$-family of eigenpairs of $H_{W,h}$ such that $\lambda_h\rightarrow E$, then the length of the zero set of $u_h$ inside the classical region $W^{-1} (]-\infty,E])$
is of order $h^{-1}$.

Although not explicitly done in \cite{LJ}, the same result is also true in the case  of $ M=\mathbb R^2$ under the condition that $\lim\inf W (x) > E\geq \inf W$, keeping the assumption that $W$ is analytic.  Let us show how we can reduce the case $M=\R^2$ to the compact situation.

\begin{proposition}
  Let us assume that $W$ is continuous and that there exists $E_1$ such that $W^{-1} (]-\infty,E_1])$ is compact.  Then the bottom of the essential spectrum of $H_{W,h}$ is bigger than $ E_1$.  Furthermore, if
   $(\lambda_h, u_h)$ is a family ($h\in ]0,h_0]$) of  eigenpairs of $H_{W,h}$ such that $\lim_{h \rightarrow 0} \lambda_h  = E_0$ with $E_0 < E_1$ and $||u_h||=1$, then given $K$ a compact neighborhood  of $W^{-1} (]-\infty,E_0])$, there exists $\epsilon_K >0$ such that
  $$
  ||u_h||_{L^2(K)} = 1  + \mathcal O \left(\exp(- \epsilon_K/h)\right)\,,
  $$
  as $h\ar 0\,$.
  \end{proposition}

This proposition is a consequence of Agmon estimates (see Helffer-Sj\"ostrand \cite{HeSj} or Helffer-Robert \cite{HeRo} for a weaker result with a remainder in  $\mathcal O_K (h^\infty)$) measuring the decay of the eigenfunctions in the classically forbidden region.  This can also be found in a weaker form in the recent book of M. Zworski \cite{Zw} (Chapter 7), which also contains
a presentation of semi-classical Carleman estimates.

Observing that in the proof of Long Jin the compact manifold $ M $ can be replaced by any compact neighborhood of $W^{-1} (]-\infty,E_0])$, we obtain:

\begin{proposition}\label{lje}
Let us assume in addition that $W$ is analytic in some compact neighborhood of $W^{-1}( ]-\infty,E_0])$, then  the length of the zero set of $u_h$ inside the classical region $W^{-1} (]-\infty,E_0])$
is of order $h^{-1}$.  More precisely, there  exist $C >0$  and $h_0 >0$ such that for all $h\in ]0,h_0]$ we have
  \begin{equation}\label{LJ-2}
  \frac 1 C  h^{-1} \leq {\rm length} \left( N(u_h) \cap W^{-1} (]-\infty,E_0]) \right) \leq C \, h^{-1}\,.
  \end{equation}
\end{proposition}
 \begin{remark}
As observed in \cite{LJ} (Remark 1.3), the results of \cite{HZZ} suggest that the behavior of the nodal sets in the classically forbidden region could be very different from the one in the classically allowed region.
\end{remark}

We can by scaling recover Proposition \ref{lbnz-P1}, and more generally treat the eigenpairs of $-\Delta_x + |x|^{2k}$.
Indeed, assume that $(-\Delta_x + |x|^{2k})u(x) = \lambda u(x)$. Write $x=\rho\,y$. Then, $(-\rho^{-2} \Delta_y + \rho^{2k} |y|^{2k} -\lambda)u(\rho y)=0$. If we choose $\rho^{2k} = \lambda$, $h = \rho^{-k-1} = \lambda^{-\frac{k+1}{2k}}$ and let $v_h(y) = h^{\frac{1}{2(k+1)}}y)\, u(h^{\frac{1}{k+1}}y)$, then, $(-h^2 \Delta_y + |y|^{2k} - 1)v_h(y) = 0$. Applying \eqref{LJ-2} to the family $v_h$ and rescaling back to the variable $x$, we find that
\begin{equation}\label{LJ-4}
\frac 1 C  \lambda^{\frac{k+2}{2k}} \leq {\rm length} \left( N(u) \cap \{x \in \R^2 ~|~ |x|^{2k} < \lambda\} \right) \leq C \, \lambda^{\frac{k+2}{2k}}\,.
\end{equation}

With this extension of Long Jin's statement, when $V= |x|^{2k}$, we also obtain an upper bound of the length of $N(u)$ in $B_V(\lambda)$.
Note that when $k\ar +\infty\,$, the problem tends to the Dirichlet problem in a ball of size $1$. We then recover that the length of $N(u)$ is of order $\sqrt{\lambda}$.

 The above method can also give results in the non-homogeneous case, at least when  \eqref{v2} is satisfied with $m_0=m_1$. We can indeed  prove the following generalization.

\begin{proposition}\label{prop8.12}~\\
Let us assume that there  exist $m\geq 1 $, $\epsilon_0 >0$ and $C >0$ such that $V$ is holomorphic in $$\mathcal D:= \{ z= (z_1,z_2) \in \mathbb C^2\,,\, | \Im z | \leq \epsilon_0 < \Re z > \}$$ and satisfies
\begin{equation}\label{lj1}
|V(z)| \leq C < \Re z >^m\,,\, \forall z \in \mathcal D\,.
\end{equation}
Suppose in addition that we have the ellipticity condition
\begin{equation}\label{lj2}
\frac 1 {C'} <x>^m \leq  V(x)\,,\, \forall x \in \mathbb R^2\,.
\end{equation}
  Then, for any $\epsilon >0$,  the length $N(u) \cap (B_V (\lambda) \setminus B_V(\epsilon \lambda))$  for an eigenpair $(u,\lambda)$ of $H_V$,  is of the order of $\lambda^{\frac 12 + \frac{2}{m}}$ as $\lambda \rightarrow + \infty$. Moreover, one can take $\epsilon =0$ when $V$ is a polynomial.
\end{proposition}

{\bf Proof}

The lower bound was already obtained by a more general direct approach in Proposition~\ref{prop7.9}. One can indeed verify  using Cauchy estimates that \eqref{lj1}  and \eqref{lj2} imply \eqref{v3} and \eqref{v2},
with  $\rho_0 = 1/2m$.
Under the previous assumptions, we consider
$$
W_\lambda (y) =\lambda^{-1} V (\lambda^{\frac 1 m} y)\,,\, v_\lambda (y) = \lambda^{\frac{1}{4m}}\, u ( \lambda^\frac 1m y)\,.
$$
We observe that  with
\begin{equation}\label{lj3a}
h = \lambda^{-\frac 12 - \frac 1m}\,,
\end{equation}
the pair $(v_\lambda,1)$ is an eigenpair for the semi-classical Schr\"odinger operator $-h^2 \Delta_y + W_\lambda (y)$:
$$
(- h^2 \Delta + W_\lambda ) v_\lambda = v_\lambda \,. $$

It remains to see  if we can extend the result of Long Jin to this situation. We essentially follow his proof, whose basic idea  goes back to Donnelly-Feffermann \cite{DF}. The difference being that $W_\lambda$ depends on $h$ through \eqref{lj3a}.\\
 The inspection of the proof\footnote{We refer here to the proof of (2.20) in \cite{LJ}.} shows that there are  three points to control.\\

 {\bf Analyticity}\\
 What we need is to have for any  $y_0$ in $\mathbb R^2\setminus \{0\}$ a complex neighborhood $\mathcal V$ of $y_0$, $h_0>0$ and  $C$ such that, for any $h\in ]0,h_0]$, $v_\lambda$ admits an holomorphic extension in $\mathcal V$
 with
 \begin{equation}\label{lj3b}
 \sup_{\mathcal V} | v_\lambda| \leq  C \exp{ \left( \frac{C}{h}\right)} \, ||v_\lambda||_{L^\infty ( \mathbb R^2)}\,.
 \end{equation}
This can be done by using the FBI transform, controlling the uniformity when $W$ is replaced by $W_\lambda$.  But this is exactly what is given by Assumption \eqref{lj1}.  Notice that this is not true in general for $y_0=0$. We cannot in general find a $\lambda$-independent neighborhhod of $0$ in $\mathbb C^2$ where $W_\lambda$ is defined and bounded. \\
  Note here
   that  $||v_\lambda||_{L^\infty ( \mathbb R^2)}$  is by standard semiclassical analysis $\mathcal O (h^{-N})$ for some $N$.\\ When $V$ is in addition a polynomial:\\
   $$
   V(x) =\sum_{j=0}^m P_j (x)
   $$
   where $P_j$ is an homogeneous polynomial of degree $j$, we get
   $$
   W_\lambda (y) = P_m (y) + \sum_{\ell =1}^{m} \lambda^{-\frac \ell m} P_{m-\ell} (y)\,,
   $$
   and we can verify the uniform analyticity property for any $y_0\,$.

 {\bf Uniform confining}\\
As we have mentioned before, Long Jin's paper was established in the case of a compact manifold  (in this case and for Laplacians, it is worth to mention the papers
of Sogge-Zelditch \cite{SZ1,SZ2})
  but it can be extended to the case of $\mathbb R^2$ under the condition that the potential is confining, the length being computed in a compact containing the classically permitted region. This is the case with $W_\lambda$. Note that if  $W_\lambda (y) \leq C_1$, then we get
  $$
  \lambda^{-1} V (\lambda^{\frac 1 m} y) \leq C_1\,,
  $$
  which implies by the ellipticity condition
  $
\frac {1}{C'}   \lambda^{-1} |\lambda^\frac 1m|^m  |y|^m \leq C_1\,,
$
that is
$$
|y| \leq (C' C_1)^{\frac 1 m}\,.
$$
{\bf Uniform doubling property}\\
Here instead of following Long Jin's proof, it is easier to refer to the results of Bakri-Casteras \cite{BaCa}, which give an explicit control in term of the $C^1$ norm of $W_\lambda$. As before, we have to use our confining assumption
 in order to establish our result in any bounded domain  $\Omega$  in $\mathbb R^2$ containing uniformly the classically permitted area $W_\lambda^{-1} (]-\infty,+1])$. This last assumption permits indeed to control the $L^2$-norm of  $v_\lambda$ from below in $\Omega$. We actually need the two following estimates (we assume \eqref{lj3a}):\\
 Given   $\Omega$ like above, for any $R>0$, there exists $C_R$ such that, for any $(x,R)$ such that $B(x, R) \subset \Omega\,$,
 \begin{equation}
 || v_\lambda||_{L^2(B(x,R))} \geq \exp{ \left( -\frac{C_R}{h}\right)} \,.
 \end{equation}
 Given   $\Omega$ like above, there exists $C$ such that, for any $(x,r)$ such that $B(x, 2r) \subset \Omega\,$,
 \begin{equation}
  || v_\lambda||_{L^2(B(x,2r))} \leq \exp{ \left( \frac{C}{h} \right)} \,  || v_\lambda||_{L^2(B(x,r))} \,.
 \end{equation}
 Here we have applied  Theorem 3.2 and Proposition 3.3 in \cite{BaCa} with electric potential $h^{-2} (W_\lambda -1)$. These two statements involve the square root of the  $C^1$ norm of the electric potential in $\overline{\Omega}$, which is $\mathcal O (h^{-1})$ in our case.\\

{\bf End of the proof}\\
Hence, considering an annulus $A(\epsilon_0, R_0)$ we get following Long Jin  that the length of the nodal set of $v_\lambda$ in this annulus is indeed of order $\mathcal O (h^{-1})$ and after rescaling
 we get the proposition for the eigenpair $(u,\lambda)$. In the polynomial case, we get the same result but in the ball $B(0,R_0)$.%
\hfill \qed\medskip

\textbf{Remarks}.\vspace{-3mm}
\begin{enumerate}
\item Long Jin's results hold in dimension $n$, not only in dimension $2$. The above extensions work in any dimension as well, replacing the length by the $(n-1)$-Hausdorff measure.
\item As observed in \cite{LJ}, the results in \cite{HZZ} suggest that the behavior of nodal sets in the classically forbidden region could be very different from the one in the classically allowed region.
\item Under the assumptions of Proposition \ref{prop8.12},  one gets from Theorem 1.1  in \cite{BaCa}  that the order of a critical point of  the zero set of an eigenfunction  of $H_V$ associated with $\lambda$ in the classically permitted region is at most of order $\lambda^{\frac 12 + \frac 1m} $.   Let us emphasize that here no assumption of analyticity for $V$  is used. On the other hand, note that using Courant's theorem and  Euler's and Weyl's formulas, one can prove that the number of critical points in the classically allowed region is at most of order  $\lambda^{1 + \frac{2}{m}}$. When $m=2$, we can verify from the results in Section \ref{S-ovals} that this upper bound cannot be improved in general.
\item For nodal sets in forbidden regions, see \cite{CaTo13}.
\end{enumerate}

\end{document}